\documentstyle[12pt]{article}

\newlength{\abstractwidth}
\renewcommand{\thefootnote}{\fnsymbol{footnote}}
\renewcommand{\thanks}[1]{\footnote{#1}} 
\newcommand{\starttext}{ \setcounter{footnote}{0}
\renewcommand{\thefootnote}{\arabic{footnote}}}

\newcommand{\be}{\begin{equation}}
\newcommand{\bea}{\begin{eqnarray}}
\newcommand{\eea}{\end{eqnarray}} 
\newcommand{\ee}{\end{equation}}
 \newcommand{\<}{\langle}
\renewcommand{\>}{\rangle}
\def\ba{\begin{eqnarray}}
\def\ea{\end{eqnarray}}

\def\o{\omega}
\def\Re{{\rm Re}}

\def\det{{\rm det}}

\def\log{\,{\rm log}\,}

\def\o{\omega}

\def\o{\omega}

\def\na{\nabla}
\def\p{\partial}

\def\ddb{{\partial\bar\partial}}

\def\na{{\nabla}}

\def\[{{\bf [}}
\def\]{{\bf ]}}

\begin{document}
\starttext \baselineskip=15pt \setcounter{footnote}{0}
\newtheorem{theorem}{Theorem}
\newtheorem{lemma}{Lemma}
\newtheorem{definition}{Definition}
\newtheorem{proposition}{Proposition}
\newtheorem{corollary}{Corollary}

\begin{center}
{\Large \bf
A FLOW OF CONFORMALLY BALANCED METRICS WITH K\"AHLER FIXED POINTS
 \footnote{Work supported in part by the National Science Foundation Grants DMS-12-66033 and DMS-1605968, and the Simons Collaboration Grant-523313. Key words: Anomaly flow, conformally balanced metric, astheno-K\"ahler metric, K\"ahler Ricci-flat metric, complex Monge-Amp\`ere flow.}}
\end{center}

\centerline{Duong H. Phong, Sebastien Picard, and Xiangwen Zhang}

\bigskip

\begin{abstract}

{\small While the Anomaly flow was originally motivated by string theory, its zero slope case is potentially of considerable interest in non-K\"ahler geometry, as it is a flow of conformally balanced metrics whose stationary points are precisely K\"ahler metrics. We establish its convergence on K\"ahler manifolds for suitable initial data. We also discuss its relation to some current problems in complex geometry.}

\end{abstract}

\section{Introduction}
\setcounter{equation}{0}

The main purpose of this paper is to study a geometric flow which is of potential interest from several viewpoints, including mathematical physics, non-K\"ahler complex geometry,  and the theory of non-linear partial differential equations. Let $X$ be a compact $n$-dimensional complex manifold, which admits a non-vanishing holomorphic $(n,0)$-form $\Omega$. The flow which we shall consider is the flow $t\to\o(t)$ of Hermitian metrics defined by
\bea
\label{AF}
\p_t(\|\Omega\|_\o\o^{n-1})
=
i\p\bar\p \o^{n-2}
\eea
with an initial data $\o_0$ which is conformally balanced, in the sense that
\bea
\label{conf_bal}
d(\|\Omega\|_{\o_0}\o_0^{n-1})=0.
\eea
Here $\|\Omega\|_{\o_0}$ is the norm of $\Omega$, defined by $\|\Omega\|_{\o_0}^2= i^{n^2} n! \, \Omega\wedge\bar\Omega \,\o_0^{-n}$. 

\smallskip
The flow (\ref{AF}) is a generalization to arbitrary dimension, with the slope parameter $\alpha'$ set to $0$, of the Anomaly flow introduced in \cite{PPZ2} for $n=3$. This flow provides a systematic approach for solving the Hull-Strominger system \cite{H1, H2, S} for supersymmetric compactifications of the heterotic string. Thus it is highly desirable to gain a better understanding for it, and the case $\alpha'=0$ provides already an important and non-trivial special case. Another defining feature of the flow (\ref{AF}) is that it preserves the conformally balanced condition (\ref{conf_bal}), and its stationary points are astheno-K\"ahler metrics (see \S 2.2 for definition). This implies that its stationary points are K\"ahler \cite{FT,MT}, so the convergence of the flow is closely related to a well-known question in non-K\"ahler geometry, namely when is a conformally balanced manifold actually K\"ahler. Also closely related is another fundamental question in non-K\"ahler geometry and algebraic-geometric stability conditions \cite{CSze, LS, Xiao1, Xiao2}, namely when does a positive $(p,p)$ cohomology class admit as representative the $p$-th power of a K\"ahler form. Finally, on K\"ahler manifolds the flow (\ref{AF}) can be viewed as a complex version of the inverse Gauss curvature flow studied extensively in convex geometry, see for example \cite{A, BIS, Ch, CT, F, Tsai}, and is in itself quite interesting as a fully non-linear partial differential equation. More details on all these motivations will be provided in Section \S 2 below.

\smallskip
In \cite{PPZ5} another generalization of the Anomaly flow to arbitrary dimension $n$ was proposed which agrees with the flow (\ref{AF}) in dimension $n=3$. For that generalization, when $\alpha'=0$ for all dimensions $n$, the flow was shown to continue to exist, as long as the curvature and torsion remain bounded and the norm of $\Omega$ does not tend to $0$. But this is all which is known at the present time. In this paper, we shall prove the following theorems:

\begin{theorem}
\label{main1} Let $X$ be an $n$-dimensional compact complex manifold with a nowhere vanishing $(n,0)$ holomorphic form $\Omega$, $n\geq 3$. Let the initial data $\omega_0$ be a Hermitian metric which is conformally balanced, i.e. the condition (\ref{conf_bal}) holds. Then the flow is parabolic, and admits in particular a unique solution on some time interval $[0,T)$ with $T>0$.
\end{theorem}

Clearly the flow (\ref{AF}) can only converge if the manifold $X$ is K\"ahler.
Conversely, if the manifold $X$ is K\"ahler, we can prove:

\begin{theorem}
\label{main2}
Consider the flow (\ref{AF}) on a compact $n$-fold $X$ equipped with a holomorphic $(n,0)$-form $\Omega$ as in the previous theorem, $n\geq 3$. Assume that the
initial data $\omega(0)$ satisfies 
\be \label{initialdata}
\| \Omega \|_{\omega(0)} \omega(0)^{n-1} = \hat{\chi}^{n-1},
\ee
where $\hat{\chi}$ is a K\"ahler metric. Then the flow (\ref{AF}) exists for all time $t>0$, and as $t \rightarrow \infty$, the solution $\omega(t)$ converges smoothly to a metric $\omega_\infty$ satisfying
\be
\omega_\infty = \| \Omega \|_{\chi_\infty}^{-2/(n-2)} \chi_\infty,
\ee
where $\chi_\infty$ is the unique K\"ahler Ricci-flat metric in the cohomology class $[\hat{\chi}]$, and $\| \Omega \|_{\chi_\infty}$ is a constant given by
\be
\| \Omega \|^2_{\chi_\infty} = {n! \over [\hat{\chi}]^n} \left( \int_X i^{n^2} \Omega \wedge \bar{\Omega} \right).
\ee
In particular, $\omega_\infty$ is K\"ahler and Ricci-flat.
\end{theorem}

As a consequence we obtain another proof of the classical theorem of Yau \cite{Y} on the existence of Ricci-flat K\"ahler metrics. A first parabolic proof of Yau's theorem was the one by Cao \cite{Cao} using the K\"ahler-Ricci flow. More recent parabolic proofs using inverse Monge-Amp\`ere flow were obtained in \cite{CK, Co1, FLM}. But, as stressed in \cite{Co1} and also discussed in Section \S 2.2 below, alternative approaches are interesting not just for alternative proofs in themselves, but also for the singularities that they would develop when no stationary point exists.

Theorem \ref{main2} will be reduced to the following theorem on a flow of complex Monge-Amp\`ere type which may be of independent interest:

\begin{theorem}
\label{main3}
Let $X$ be a compact K\"ahler manifold with K\"ahler metric $\hat{\chi}$. Let $f \in C^\infty(X,{\bf R})$ be a given function. Consider the flow
\be \label{MA-flow}
\p_t \varphi = e^{-f} \, {(\hat{\chi} + i \ddb \varphi)^n \over \hat{\chi}^n}, \ \ \varphi(x,0)=0,
\ee
where $\varphi(t)$ is subject to the plurisubharmonicity condition
\be \label{positivity}
\hat{\chi} + i \ddb \varphi >0.
\ee
Then the flow $\varphi(t)$ exists for all $t>0$, and the averages $\tilde\varphi(t)=\varphi(t)-{1 \over V} \int_X \varphi\hat\chi^n$, with $V = \int_X \hat{\chi}^n$, converge in $C^\infty$.
\end{theorem}

The flow (\ref{MA-flow}) shares with the K\"ahler-Ricci flow and the inverse Monge-Amp\`ere flow ($MA^{-1}$-flow) introduced by Collins-Hisamoto-Takahashi \cite{Co1} the fact that the Monge-Amp\`ere determinant ${(\hat{\chi} + i \ddb \varphi)^n /\hat{\chi}^n}$ appears in the right hand side. However, it appears here as a first power, and not as a log or as an inverse power, which makes it not concave. In particular, it does not fall within the scope of  standard parabolic PDE methods such as those developed in \cite{PT}.

\medskip

The paper is organized as follows. In \S 2, we provide some background and motivations for the flow (\ref{AF}). In \S 3, we compute the evolution equation for the metric and prove the short time existence as claimed in Theorem \ref{main1}. In \S 4, the Anomaly flow is reduced to the study of a complex Monge-Amp\`ere flow for initial data satisfying (\ref{initialdata}). We establish all the necessary estimates for the long-time behavior of the flows. In \S 5, we prove the convergence of the flow which completes the proof of Theorems \ref{main2} and \ref{main3}. In \S 6, we provide some further remarks about the Anomaly flow and its possible connection to some interesting problems in complex geometry. The appendix contains some calculations used in \S 3.

\section{Motivations for the flow}
\setcounter{equation}{0}

We  provide now some details on the three different contexts which make the flow (\ref{AF}) of particular interest.

\subsection{The Hull-Strominger system}

Let $X$ be a compact $3$-fold equipped with a nowhere vanishing holomorphic $(3,0)$-form $\Omega$. Let $t\to \Phi(t)$ be a flow of $(2,2)$-forms with $d\Phi=0$ for any $t$. Then the following flow was introduced in \cite{PPZ2}
\bea
\label{AF-o}
\p_t(\|\Omega\|_\o\o^2)
=
i\p\bar\p\o-{\alpha'\over 4}({\rm Tr}(Rm \wedge Rm)-\Phi(t))
\eea
with initial data $\o_0$ required to satisfy the condition $d(\|\Omega\|_{\o_0}\o_0^2)=0$. Here $\alpha'$ is a physical parameter called the slope, and $Rm$ is the Riemann curvature of the metric $\o$, viewed as a $(1,1)$-form valued in the bundle of endomorphisms of $T^{1,0}(X)$. One case of particular interest is when 
$X$ is also equipped with a holomorphic vector bundle $E\to X$ with $c_1(E)=0$,
and $\Phi(t)={\rm Tr}(F\wedge F)$ where $F$ is the curvature of a Hermitian metric $H(t)$ on $E$ which itself evolves simultaneously under the Donaldson heat flow
\bea
\label{AF-H}
H^{-1}\p_tH=-(\Lambda_{\omega} iF), \qquad H(0)=H_0.
\eea
Here $H_0$ is a given initial metric, $F$ is viewed as a $(1,1)$-form valued in the bundle of endomorphisms of $E$, and $(\Lambda_{\omega} iF)^\alpha{}_\beta=g^{j\bar k}F_{\bar kj}{}^\alpha{}_\beta$ is the usual Hodge contraction. The stationary points of the simultaneous flows (\ref{AF-o}) and $(\ref{AF-H})$ are precisely the solutions of a system of equations proposed independently by C. Hull \cite{H1,H2} and A. Strominger \cite{S} for supersymmetric compactifications of the heterotic string to $M^{3,1}\times X$, where $M^{3,1}$ is Minkowski space-time. The Hull-Strominger system of equations is a generalization of the more specific case considered earlier by Candelas, Horowitz, Strominger and Witten \cite{CHSW}, where $E$ was set to $T^{1,0}(X)$. If we consider only the stationary points of the system, we can set $H=\o$, and consequently $\Phi={\rm Tr}(Rm\wedge Rm)$. In this case the term ${\rm Tr}(Rm \wedge Rm)-\Phi(t)$ vanishes trivially, and it was shown in \cite{CHSW} that the stationary points of both (\ref{AF-o}) and (\ref{AF-H}) are given by the same condition of $\o$ being a Ricci-flat K\"ahler metric.

\medskip
A major reason for considering the flow (\ref{AF-o}), among all the flows which have the same stationary points, is that it preserves the condition that the form $\|\Omega\|_\o\o^2$ be closed. In fact, the $(2,2)$ de Rham class of $\|\Omega\|_{\o(t)}\o(t)^2$ can be easily read off from the flow. This is however only a weak substitute for the analogous statement for the $(1,1)$ K\"ahler class of $\o$ in the K\"ahler-Ricci flow, and this accounts for many new difficulties that one encounters in the study of (\ref{AF-o}).

\smallskip
The flow (\ref{AF-o}) can be expressed more explicitly in terms of a flow for the $(1,1)$-form $\o(t)$ itself \cite{PPZ5}. Perhaps surprisingly, this formulation shows that (\ref{AF-o}) can be viewed as a generalization of the Ricci flow with a triple complication, namely the metrics are not K\"ahler (or Levi-Civita), the norms $\|\Omega\|_\o$ also occur, and so do quadratic expressions in the curvature tensor.

\smallskip
The flow (\ref{AF-o}) was shown in \cite{PPZ6} to produce a simpler and unified way of recovering the solutions of the Hull-Strominger system on Calabi-Eckmann-Goldstein-Prokushkin fibrations \cite{CE, GP} originally found by Fu and Yau \cite{FY1, FY2}. It also inspired directly the solutions of the Fu-Yau Hessian equations found recently in \cite{PPZ8}, which extended our earlier work \cite{PPZ1, PPZ4}. However, it is still largely unexplored, and its long-term behavior is already known to exhibit in general a rather intricate behavior, such as a particular sensitivity to the initial data \cite{FHP, PPZ6, PPZ7}.

\smallskip
Because of all these complications, it is necessary to examine the flow (\ref{AF-o}) in some simpler settings. 
The above choice of $\Phi$ in \cite{CHSW} has the effect of eliminating the terms in the right hand side of (\ref{AF-o}) which are quadratic in the curvature tensors. This effect is also achieved simply by considering the case $\alpha'=0$. The flow (\ref{AF-o}) can then be considered on its own. And while the dimension $3$ for $X$ is required for the interpretation of the Hull-Strominger system as a compactification of the heterotic string, from the pure mathematical standpoint, we can consider all dimensions $n\geq 3$, and the flow we propose in (\ref{AF}) is the natural generalization of the Anomaly flow (\ref{AF-o}) with $\alpha'=0$.

\subsection{Generalizations of the K\"ahler condition}
\setcounter{equation}{0}

The flow (\ref{AF}) also fits into the broad question of when a compact complex manifold $X$ may admit a K\"ahler metric or a weaker substitute. We review some of the relevant notions.

\smallskip

Let $(X, J)$ be a compact complex manifold of complex dimension $n$ equipped with a Hermitian metric $g$, and denote its K\"ahler form by $\omega=i\sum g_{\bar k j} dz^j \wedge d\bar z^k$. Since the K\"ahler form of a Hermitian metric determines the metric (as $J$ is fixed), by an abuse of terminology we will not distinguish between the two notions.

If $d\omega=0$, then $g$ is called a K\"ahler metric, and $X$ is called a K\"ahler manifold. For each $2\leq k\leq n-1$, it 
is natural to consider a weaker condition of the form
\be
d\, \omega^{k} =0.
\ee
The case $k=n-1$ was introduced by Michelsohn \cite{Michelsohn} who called such a metric balanced, and the manifold $X$ a balanced manifold. By an observation due to Gray and Hervella \cite{GH} conditions such as $d\omega^k=0$ for some $2\leq k \leq n-2$ actually imply that $\o$ is K\"ahler, so the case $k=n-1$ is the only non-trivial generalization of the K\"ahler property. Michelsohn found an intrinsic characterization of compact manifolds with balanced metrics by means of positive currents. Using such a characterization, Alessandrini and Bassanelli proved that the existence of balanced metrics is preserved under birational transformations in \cite{AB}. Remarkably, as we just saw in the previous section \S 2.1, balanced metrics also arise in string theory since the torsion constraint equation $d(\|\Omega\|_\o\o^2)=0$ in dimension $3$ of the Hull-Strominger system just means that the metric $\|\Omega\|_\o^{1\over 2}\omega$ is balanced, see \cite{LY}. The existence of balanced metrics on compact Hermitian manifolds has been studied extensively (see e.g. \cite{Fei1,Fei2,FIUV,FG,Fu,FLY,FWW1,Michelsohn, TW, Ug} and references therein).

\medskip

On the other hand, another natural generalization of the K\"ahler condition is
\be
i\partial\bar\partial \omega^\ell =0, 
\ee
for some $\ell$ with $1\leq\ell \leq n-1$. If $\ell=n-1$, then $\omega$ is called a Gauduchon metric. In \cite{Gauduchon}, Gauduchon proved that there always exists a unique Gauduchon metric, up to a constant conformal factor, in the conformal class of a given Hermitian metric. Gauduchon manifolds also provide a natural setting for an extension \cite{LY0, LT} of the Donaldson-Uhlenbeck-Yau theory of Hermitian-Yang-Mills metrics on stable bundles over K\"ahler manifolds. 

\smallskip

If $\ell=n-2$, then $\omega$ is called an astheno-K\"ahler metric. This notion was introduced by Jost and Yau in \cite{JY} to establish the existence of Hermitian harmonic maps, and it turns out to be particularly interesting for many analytic arguments to be useful. For example, Tosatti and Weinkove proved in \cite{TW1} the Calabi-Yau theorems for Gauduchon and strongly Gauduchon metrics on the class of compact astheno-K\"ahler manifolds. Recently, the existence such metrics on compact complex manifolds has been studied widely, see for example \cite{FGV, FT, LU, Matsuo, MT}.

\smallskip

An interesting question raised in \cite{STW} asks whether a compact complex non-K\"ahler manifold can admit both an astheno-K\"ahler metric and a balanced metric. Very recently, such examples were constructed in \cite{FGV,LU}. In fact, the astheno-K\"ahler metric of Latorre-Ugarte \cite{LU} is $k$-th Gauduchon \cite{FWW2}, meaning that $i \ddb \omega^k \wedge \omega^{n-k-1} =0$, for every $1 \leq k \leq n-1$. However, a single metric on a compact manifold cannot be both balanced and astheno-K\"ahler, unless it is K\"ahler. This statement can be viewed as a generalization to arbitrary dimensions of the arguments of \cite{CHSW}. It was proved in Matsuo-Takahashi \cite{MT}, and extended in Fino-Tomassini \cite{FT} to the case of metrics which are both conformally balanced and astheno-K\"ahler. For easy reference, we state these results as a lemma, and include an alternative proof in this paper (see Corollary \ref{stationary-pts} in \S 3 for (ii) and (e) in \S 6 for (i)):

\begin{lemma} \cite{FT, MT}
\label{cf+ak}
Let $X$ be an $n$-dimensional compact complex manifold with Hermitian metric $\o$. Then $\o$ is K\"ahler if one of the following conditions is satisfied

{\rm (i)} $\o$ satisfies both $d\o^{n-1}=0$ and $i\p\bar\p\o^{n-2}=0$,

{\rm (ii)} $X$ admits a nowhere vanishing holomorphic $(n,0)$-form $\Omega$ and $\o$ satisfies both $d(\|\Omega\|_\o\o^{n-1})=0$ and $i\p\bar\p\o^{n-2}=0$. In this case, $\o$ is also Ricci-flat.

\end{lemma}

We return now to the flow (\ref{AF}). Since the right hand side of the flow is both $d$ and $i\p\bar\p$ exact, it follows that if $\|\Omega\|_{\o_0}\o_0^{n-1}$ is $d$-closed or $i\p\bar\p$-closed, then $\|\Omega\|_{\o}\o^{n-1}$ remains $d$-closed or $i\p\bar\p$-closed along the flow, and the de Rham cohomology class $[\|\Omega\|_{\o_0}\o_0^{n-1}]$ or the Bott-Chern cohomology class $[\|\Omega\|_{\o_0}\o_0^{n-1}]_{BC}$ is preserved.

\smallskip
Moreover, thanks to (ii) in Lemma \ref{cf+ak}, we see that
if the Anomaly flow (\ref{AF}) converges, the limit metric must be a K\"ahler Ricci-flat metric. Therefore, the flow provides a deformation path in the space of conformally balanced metrics to a K\"ahler metric. If no K\"ahler metric exists on $X$, then the flow cannot converge, and either its singularities in finite-time or long-term behavior should provide an analytic measure of the absence of K\"ahler metrics.

\subsection{Parabolic fully non-linear equations on Hermitian manifolds}

The theory of parabolic fully non-linear equations on Hermitian manifolds has been developed extensively over the years, and has resulted in some powerful and general results. As an example, let $(X,\o)$ be a compact Hermitian manifold equipped with a smooth closed form $\hat\chi$, and consider the equation
\bea
\label{AF-u}
\p_t u=F(A(i\p\bar\p u))-\psi(z),
\qquad
u(0)=u_0. 
\eea
Here $i\p\bar\p u$ is the complex Hessian of $u$, and $A(i\p\bar\p u)$ is the set of eigenvalues of the form $\hat\chi+i\p\bar\p u$ with respect to the metric $\o$. The function $F(\lambda)$ is a given function defined on a given cone $\Gamma$, and $A(i\p\bar\p u)$ is required to be in $\Gamma$ for all time. 

\smallskip

The function $F$ is required to satisfy many conditions, including conditions amounting to the parabolicity of the flow (\ref{AF-u}). But an additional condition, originating from the earliest works of Caffarelli-Nirenberg-Spruck \cite{CNSIII} and Krylov \cite{Kr0} on the theory of fully non-linear equations, is that $F$ be a concave function of $\lambda\in \Gamma$. This concavity condition, together with the existence of subsolutions and more technical hypotheses, is now known to lead to some very general theorems which can apply to a wide variety of geometric flows that have been studied in the literature (see e.g. \cite{PT} and references therein).

\smallskip
The Anomaly flow on Calabi-Eckmann-Goldstein-Prokushkin fibrations studied in \cite{PPZ6} provides however a natural geometric example of a flow that does not satisfy the concavity condition, but which is nevertheless well-behaved. This suggests the existence of a useful theory going beyond concave equations, the formulation of which would require the treatment of many more examples. As we shall see, the flow (\ref{AF}), in its realization (\ref{MA-flow}), provides another instructive example. Note that it corresponds to $F=e^{f(z)} \prod_{j=1}^n\lambda_j$, unlike the K\"ahler-Ricci flow, which corresponds to $F=\log \prod_{j=1}^n\lambda_j$ which is concave.

\section{Short-time existence and proof of Theorem \ref{main1}}
\setcounter{equation}{0}

To prove Theorem \ref{main1}, 
we need a lemma which will allow us to compute relevant quantities with the Hodge star operator. 
For a $(p,q)$-form $\Theta$ on a manifold $X$, its components $\Theta_{\bar k_1\cdots\bar k_q j_1\cdots j_p}$ are defined by
\bea
\Theta=
{1\over p!q!}
\sum \Theta_{\bar k_1\cdots\bar k_q j_1\cdots j_p}\,
dz^{j_p}\wedge\cdots\wedge dz^{j_1}\wedge
d\bar z^{k_q}\wedge\cdots\wedge d\bar z^{k_1}.
\eea
For a $(p,p)$-form $\Theta$, we define
\bea\label{defineTr}
{\rm Tr}\,\Theta=\<\Theta, \o^p\>=
i^{-p}\, \prod_{\ell=1}^p g^{k_\ell\bar j_\ell}\,\Theta_{\bar j_1k_1\cdots\bar j_pk_p}.
\eea

\begin{lemma} \label{star-wedges}
Let $\alpha \in \Omega^{1,1}(X,{\bf R})$, $\Phi \in \Omega^{2,2}(X,{\bf R})$ and $\Psi \in \Omega^{3,3}(X,{\bf R})$ on a complex manifold $X$ of dimension $n \geq 3$. Let $\omega = i g_{\bar{k} j} dz^j \wedge d \bar{z}^k$ be a Hermitian metric and $\star$ its associated Hodge star operator. Then
\be \label{star-wedge-n-2}
\star ( \alpha \wedge \omega^{n-2}) = - (n-2)! \, \alpha + (n-2)! ({\rm Tr} \, \alpha) \, \omega,
\ee
\be \label{star-wedge-n-3}
(\star (\Phi \wedge \omega^{n-3}) )_{\bar{k} j} = i(n-3)! g^{s \bar{r}} \Phi_{\bar{r} s \bar{k} j} + i{(n-3)! \over 2}  ( {\rm Tr} \, \Phi ) \, g_{\bar{k} j},
\ee
and if $n \geq 4$, then
\be \label{star-wedge-n-4}
(\star (\Psi \wedge \omega^{n-4}) )_{\bar{k} j} = {(n-4)! \over 2} g^{q \bar{p}} g^{s \bar{r}} \Psi_{\bar{r} s \bar{p} q \bar{k} j} + i {(n-4)! \over 6}  ( {\rm Tr} \, \Psi ) \, g_{\bar{k} j}.
\ee
\end{lemma}

{\it Proof:} We will use the general formula for the Hodge star operator on $(n-1,n-1)$ forms, as given in (Lemma 3, \cite{PPZ5}). For any $(n-1,n-1)$ form $\Theta$, there holds
\bea
(\star\Theta)_{\bar j k}
=
{1\over (n-1)!(n^2-6n+6)}
\bigg\{
6\, i^{-(n-2)}\, \prod_{p=1}^{n-2}g^{k_p\bar j_p}\Theta_{\bar jk\bar j_1k_1\cdots\bar j_{n-2}k_{n-2}}
+(n-6)({\rm Tr}\Theta)\, i \, g_{\bar j k}\bigg\}.
\nonumber
\eea
In fact, the proof of (\ref{star-wedge-n-2}) can also be found in \cite{PPZ5, Huy}, but we provide details here for completeness. First, we notice that
\be
(\alpha \wedge \omega^{n-2})_{\bar{k} j \bar{k_1} j_1 \cdots \bar{k}_{n-2} j_{n-2}} = i^{n-2} \alpha_{ \{ \bar{k} j} g_{\bar{k_1} j_1} \cdots g_{\bar{k}_{n-2} j_{n-2} \}},
\ee
\be
(\Phi \wedge \omega^{n-3})_{\bar{k} j \bar{r} s \bar{k_1} j_1 \cdots \bar{k}_{n-3} j_{n-3}} = i^{n-3} \Phi_{ \{ \bar{k} j \bar{r} s} g_{\bar{k_1} j_1} \cdots g_{\bar{k}_{n-3} j_{n-3} \}},
\ee
\be
(\Psi \wedge \omega^{n-4})_{\bar{k} j \bar{r} s \bar{p} q \bar{k_1} j_1 \cdots \bar{k}_{n-4} j_{n-4}} = i^{n-4} \Psi_{ \{ \bar{k} j \bar{r} s \bar{p} q} g_{\bar{k_1} j_1} \cdots g_{\bar{k}_{n-4} j_{n-4} \}}.
\ee
Here we use the notation $\{ \}$ to mean antisymmetrization in both the barred and unbarred indices. Using the formula for the Hodge star on $(n-1,n-1)$ forms exhibited above, we deduce that we must have
\be \label{star-n-2-unknown}
(\star (\alpha \wedge \omega^{n-2}) )_{\bar{k} j} = a_n \alpha_{\bar{k} j} + i b_n ({\rm Tr} \, \alpha) g_{\bar{k} j},
\ee
\be \label{star-n-3-unknown}
(\star (\Phi \wedge \omega^{n-3}) )_{\bar{k} j} = i c_n g^{s \bar{r}} \Phi_{\bar{r} s \bar{k} j} + i d_n ({\rm Tr} \, \Phi) g_{\bar{k} j},
\ee
\be \label{star-n-4-unknown}
(\star (\Psi \wedge \omega^{n-4}) )_{\bar{k} j} = e_n g^{q \bar{p}} g^{s \bar{r}} \Psi_{\bar{r} s \bar{p} q \bar{k} j} + i f_n ({\rm Tr} \, \Psi) g_{\bar{k} j},
\ee
for some coefficients $a_n, b_n, c_n, d_n, e_n, f_n$ to be determined. If ${\rm Tr} \, \alpha = 0$, a direct computation shows
\be
\star (\alpha \wedge \omega^{n-2}) = - (n-2)! \alpha.
\ee
Therefore $a_n = -(n-2)!$. Taking $\alpha = \omega$ and using $\star \omega^{n-1} = (n-1)! \omega$, we deduce the relation
\be
(n-1)! = a_n + n b_n.
\ee
It follows that $b_n = (n-2)!$ and we have established (\ref{star-wedge-n-2}). 

Next, we test $\Phi = \alpha \wedge \omega$, which has components
\be
(\alpha \wedge \omega)_{\bar{r} s \bar{k} j} = i (\alpha_{\bar{k} j} g_{\bar{r} s} - \alpha_{\bar{k} s} g_{\bar{r} j} - \alpha_{\bar{r} j} g_{\bar{k} s} + \alpha_{\bar{r} s} g_{\bar{k} j}). 
\ee
Thus
\be
g^{s \bar{r}} (\alpha \wedge \omega)_{\bar{r} s \bar{k} j} = (n-2) i \alpha_{\bar{k} j} - ({\rm Tr} \, \alpha) g_{\bar{k} j}.
\ee
and
\be
{\rm Tr} \, (\alpha \wedge \omega) = 2 (n-1) ({\rm Tr} \, \alpha).
\ee
Equating (\ref{star-n-2-unknown}) with (\ref{star-n-3-unknown}), we therefore we have the relation
\bea
-(n-2)! \alpha_{\bar{k} j} + i (n-2)! ({\rm Tr} \, \alpha) g_{\bar{k} j}= i^2 c_n (n-2) \alpha_{\bar{k} j} - i c_n ({\rm Tr} \, \alpha) g_{\bar{k} j} + 2 i d_n (n-1) ({\rm Tr} \, \alpha) g_{\bar{k} j}.\nonumber
\eea
Taking $\alpha$ with ${\rm Tr} \, \alpha = 0$, we see that $c_n = (n-3)!$. Solving for $d_n$ gives $d_n = {(n-3)! \over 2}$.
This establishes (\ref{star-wedge-n-3}). 

We now test $\Psi = \Phi \wedge \omega$. Its components can be worked out to be
\bea
(\Phi \wedge \omega)_{\bar{r} s \bar{p} q \bar{k} j} &=& i  \bigg\{ \Phi_{\bar{r} s \bar{p} q} g_{\bar{k} j} + \Phi_{\bar{r} j \bar{p} s} g_{\bar{k} q} + \Phi_{\bar{r} q \bar{p} j} g_{\bar{k} s}+ \Phi_{\bar{k} s \bar{r} q} g_{\bar{p} j} + \Phi_{\bar{k} j \bar{r} s} g_{\bar{p} q} \nonumber\\
&& + \Phi_{\bar{k} q \bar{r} j} g_{\bar{p} s} + \Phi_{\bar{p} s \bar{k} q} g_{\bar{r} j} + \Phi_{\bar{p} j \bar{k} s} g_{\bar{r} q} + \Phi_{\bar{p} q \bar{k} j} g_{\bar{r} s} \bigg\}. 
\eea
Then
\be
g^{q \bar{p}} g^{s \bar{r}} (\Phi \wedge \omega)_{\bar{r} s \bar{p} q \bar{k} j} =  2 i (n-3)  g^{s \bar{r}} \Phi_{\bar{k} j \bar{r} s}- i {\rm Tr} \, \Phi \, g_{\bar{k} j},
\ee
and
\be
{\rm Tr} \, (\Phi \wedge \omega) = 3(n-2) ({\rm Tr} \, \Phi).
\ee
Equating (\ref{star-n-3-unknown}) and (\ref{star-n-4-unknown}) with $\Psi = \Phi \wedge \omega$ gives the relation
\bea \label{e_n-f_n}
& \ & i(n-3)! g^{s \bar{r}} \Phi_{\bar{r} s \bar{k} j} + i{(n-3)! \over 2}  ( {\rm Tr} \, \Phi ) \, g_{\bar{k} j} \nonumber\\
&=&  2 i (n-3)  e_n g^{s \bar{r}} \Phi_{\bar{k} j \bar{r} s} - i e_n( {\rm Tr} \, \Phi) \, g_{\bar{k} j} + 3 i (n-2) f_n ({\rm Tr} \, \Phi) g_{\bar{k} j}.
\eea
If we choose $\Phi$ with ${\rm Tr} \, \Phi = 0$, we see that $e_n = {(n-4)! \over 2}$. Substituting the value of $e_n$ into (\ref{e_n-f_n}), we deduce $f_n = {(n-4)! \over 6}$.
Q.E.D.

\

\par 
Before stating the full evolution of the components of the metric in the Anomaly flow (\ref{AF}), we establish notation which will be used subsequently. 
 The torsion of $\omega(t)$ is given by $T = i \p \omega$ and $\bar{T}= - i \bar{\p} \omega$, with components
\bea
T={1\over 2} T_{\bar{k} j \ell} \, dz^\ell \wedge dz^j \wedge d\bar z^k,
\quad
\bar T={1\over 2} \bar{T}_{k \bar{j} \bar{\ell}} \, d\bar z^\ell\wedge d\bar z^j\wedge dz^k,
\eea
given by
\bea
T_{\bar{k} j \ell}=\p_jg_{\bar{k} \ell}-\p_\ell g_{\bar{k} j},
\quad
\bar T_{k \bar{j} \bar{\ell}}=\p_{\bar{j}}g_{\bar{\ell} k}-\p_{\bar{\ell}}g_{\bar{j} k}.
\eea
The torsion $1$-form $\tau$ is given by
\be
\tau = T_\ell dz^\ell, \ \ \ T_\ell = g^{j \bar{k}} T_{\bar{k} j \ell}.
\ee
We may take the norms of $T$ and $\tau$ by setting
\be
|T|^2 = g^{m \bar{k}} g^{j \bar{n}} g^{\ell \bar{p}} T_{\bar{k} j \ell} \bar{T}_{m \bar{n} \bar{p}}, \ \ \ |\tau|^2 = g^{j \bar{k}} T_j \bar{T}_{\bar{k}}.
\ee
The curvature tensor of the Chern connection $\na$ of $\omega$ is given by
\bea
R_{\bar kj}{}^p{}_q=-\p_{\bar k}(g^{p\bar \ell}\p_jg_{\bar\ell q}),
\
Rm = R_{\bar kj}{}^p{}_q \, dz^j \wedge d\bar z^k
\in \Omega^{1,1}(X)\otimes End(T^{1,0}(X)).
\eea
For a general Hermitian metric $\omega$, there are 4 notions of Ricci curvature.
\bea\label{def-Ric}
R_{\bar kj}=R_{\bar kj}{}^p{}_p,\qquad \tilde R_{\bar kj}=R^p{}_p{}_{\bar kj},
\qquad
R_{\bar kj}'=R_{\bar k}{}^p{}_p{}_j,
\qquad
R_{\bar kj}''=R^p{}_{j\bar k}{}_p.
\eea
Conformally balanced metrics have additional structure, and their curvatures satisfy (see Lemma 5 in \cite{PPZ5})
\be \label{conf-bal-riccis}
R_{\bar kj}'=R_{\bar kj}''= {1 \over 2} R_{\bar kj}, \ \ \tilde R_{\bar kj}= {1 \over 2} R_{\bar kj}+\na^mT_{\bar kjm}.
\ee

\

We can now state the formula for the evolution of the metric tensor along the Anomaly flow (\ref{AF}). The structure of the torsion terms can be compared to other geometric flows studied in different settings e.g. \cite{BV, Bry, Ka, LoWe, ST,ST2, TW2}. The appearance of quadratic torsion terms proportional to $g_{\bar{k} j}$ when $n \geq 4$ is reminiscent of evolution of the metric in the Laplacian flow for closed $G_2$ structures \cite{Bry,Ka,LoWe}. 
\par The difference in expressions when $n=3$ and $n \geq 4$ is due to the appearance of $(n-2)(n-3) i \p \omega \wedge \bar{\p} \omega \wedge \omega^{n-3}$ when expanding $i \ddb \omega^{n-2}$ if $n \geq 4$. In deriving (\ref{metric-evol-n>3}), we cancelled coefficients of the form ${(n-3) \over (n-3)}$, making it invalid to substitute $n=3$ in this expression.

\begin{theorem}
\label{main4}
Start the Anomaly flow ${d \over dt} (\| \Omega \|_\omega \omega^{n-1} )= i \ddb (\omega^{n-2})$ with initial metric $\omega_0$ satisfying $d (\| \Omega \|_{\omega_0} \omega_0^{n-1} )= 0$. If $n=3$, the evolution of the metric is given by
\bea
\p_t g_{\bar kj}
 =
{1\over 2\|\Omega\|_\o}\bigg[-\tilde R_{\bar kj} + g^{m \bar{\ell}} g^{s \bar{r}} T_{\bar{r} m j} \bar{T}_{s \bar{\ell} \bar{k}} \bigg].
\eea
If $n \geq 4$, the evolution of the metric is given by
\bea \label{metric-evol-n>3}
\p_t g_{\bar{k} j} &=& {1 \over (n-1) \| \Omega \|_\omega} \bigg[ -\tilde{R}_{\bar{k} j} + {1 \over 2 (n-2)} (|T|^2- 2 |\tau|^2) \, g_{\bar{k} j} \nonumber\\
&&- {1 \over 2} g^{q \bar{p}} g^{s \bar{r}} T_{\bar k q s} \bar{T}_{j\bar p\bar r}
+ g^{s\bar r} (T_{\bar k j s} \bar{T}_{\bar r} + T_s \bar{T}_{j\bar k \bar r})  +  T_j \bar{T}_{\bar k} \bigg].
\eea
\end{theorem}

As a immediate consequence, we see that the flow has a short-time solution. Indeed, by definition $\tilde{R}_{\bar{k} j} = - g^{p \bar{q}} \p_{\bar{q}} \p_p g_{\bar{k} j} + g^{p \bar{q}} g^{r \bar{s}} \p_{\bar{q}} g_{\bar{k} r} \p_p g_{\bar{s} j}$. The leading term of the linearization is the different operator
\bea
\delta g_{\bar kj}
\to 
-g^{p\bar q}\p_p\p_{\bar q}\delta g_{\bar kj}
\eea
which is elliptic. Hence the flow is strictly parabolic.  Theorem \ref{main1} is proved.

\bigskip

\par 
\noindent {\it Proof:} The evolution equation when $n=3$ was already given in \cite{PPZ5}, so we assume that $n \geq 4$. Since $\| \Omega \|_\omega^2 = \Omega \bar{\Omega} (\det g)^{-1}$, we have
\be
{d \over dt} \| \Omega \|_\omega = - {1 \over 2} \| \Omega \|_\omega {\rm Tr} \, \dot{\omega}.
\ee
The Anomaly flow (\ref{AF}) can be written as
\bea \label{af-expand}
& \ & - {1 \over 2} \| \Omega \|_\omega ({\rm Tr} \, \dot{\omega}) \omega^{n-1} + (n-1) \| \Omega \|_\omega \dot{\omega} \wedge \omega^{n-2} \nonumber\\
&=& (n-2) i \ddb \omega \wedge \omega^{n-3} + i (n-2)(n-3) T \wedge \bar{T} \wedge \omega^{n-4}.
\eea
Wedging both sides by $\omega$ gives the following equation of top forms
\bea
& \ & - {1 \over 2} \| \Omega \|_\omega ({\rm Tr} \, \dot{\omega}) \omega^{n} + (n-1) \| \Omega \|_\omega \dot{\omega} \wedge \omega^{n-1} \nonumber\\
&=& (n-2) i \ddb \omega \wedge \omega^{n-2} + i (n-2)(n-3) T \wedge \bar{T} \wedge \omega^{n-3}.
\eea
Using the identities (\ref{1-1-contract}), (\ref{2-2-contract}), (\ref{3-3-contract}) for contracting forms, we obtain
\bea
& \ & - {1 \over 2} \| \Omega \|_\omega ({\rm Tr} \, \dot{\omega}) + {(n-1) \over n}  \| \Omega \|_\omega ({\rm Tr} \, \dot{\omega}) \nonumber\\
& = & - {(n-2) \over 2n(n-1)} g^{j \bar{k}} g^{\ell \bar{m}} (i \ddb \omega)_{\bar{k} j \bar{m} \ell}  - {(n-3)\over 6n(n-1)}  g^{j \bar{k}} g^{q \bar{p}} g^{s \bar{r}} (T \wedge \bar{T})_{\bar{k} j \bar{p} q \bar{r} s}.
\eea
Therefore
\be \label{trace-iddb-omega}
({\rm Tr} \, \dot{\omega}) = {1 \over (n-1) \| \Omega \|_\omega}{\rm Tr} \, (i \ddb \omega) + {i (n-3) \over 3 (n-1)(n-2)} {1 \over \| \Omega \|_\omega} {\rm Tr} \, (T \wedge \bar{T}).
\ee
We now apply the Hodge star operator with respect to $\omega$ to (\ref{af-expand})
\bea
& \ & - {1 \over 2} \| \Omega \|_\omega ({\rm Tr} \, \dot{\omega})  \star \omega^{n-1} + (n-1) \| \Omega \|_\omega  \star (\dot{\omega} \wedge \omega^{n-2}) \nonumber\\
&=& (n-2) \star (i \ddb \omega \wedge \omega^{n-3}) + i (n-2)(n-3) \star (T \wedge \bar{T} \wedge \omega^{n-4}).
\eea
Substituting Lemma \ref{star-wedges},
\bea
& \ & - {(n-1) \over 2} \| \Omega \|_\omega ({\rm Tr} \, \dot{\omega}) g_{\bar{k} j} + (n-1) \| \Omega \|_\omega \{ -  \p_t g_{\bar{k} j} + ({\rm Tr} \, \dot{\omega})  g_{\bar{k} j} \}\\
&=&  g^{s \bar{r}} (i \ddb \omega)_{\bar{r} s \bar{k} j} +  {1 \over 2} ({\rm Tr} \, i \ddb \omega ) g_{\bar{k} j} +  {1 \over 2} g^{q \bar{p}} g^{s \bar{r}} (T \wedge \bar{T})_{\bar{r} s \bar{p} q \bar{k} j}   + i {1 \over 6} {\rm Tr} \, (T \wedge \bar{T}) g_{\bar{k} j}.   \nonumber
\eea
Substituting (\ref{trace-iddb-omega}), we see that the $({\rm Tr} \, i \ddb \omega )$ terms cancel and we are left with
\bea \label{p-t-metric1}
\p_t g_{\bar{k} j} &=& - {1 \over (n-1) \| \Omega \|_\omega} g^{s \bar{r}} (i \ddb \omega)_{\bar{r} s \bar{k} j}  - {1 \over 2 (n-1) \| \Omega \|_\omega} g^{q \bar{p}} g^{s \bar{r}} (T \wedge \bar{T})_{\bar{r} s \bar{p} q \bar{k} j} \nonumber\\
&&- {i \over 6 (n-1)(n-2)\| \Omega \|_\omega } {\rm Tr} \, (T \wedge \bar{T} ) \, g_{\bar{k} j}.
\eea

Next, we compute the components of the terms of the right-hand side. By substituting the relation of Ricci curvatures for conformally balanced metrics (\ref{conf-bal-riccis}) into identity (\ref{tr-iddb-omega}) for components of $i \ddb \omega$ from the appendix, we obtain
\be 
g^{s \bar{r}} (i \ddb \omega)_{\bar{r} s \bar{k} j} = \tilde{R}_{\bar{k} j} - g^{m \bar{\ell}} g^{s \bar{r}} T_{\bar{r} m j} \bar{T}_{s \bar{\ell} \bar{k}}. 
\ee
Substituting this identity together with (\ref{tr-tr-TT}) and (\ref{tr-tr-tr-TT}) into (\ref{p-t-metric1}) gives equation (\ref{metric-evol-n>3}) for the evolution of the metric. Q.E.D.

\

As a consequence of the above calculation, we obtain a proof for case (ii) in Lemma \ref{cf+ak}.

\begin{corollary} \label{stationary-pts}
Stationary points of the Anomaly flow ${d \over dt} (\| \Omega \|_\omega \omega^{n-1} )= i \ddb (\omega^{n-2})$ with initial metric $\omega_0$ satisfying $d (\| \Omega \|_{\omega_0} \omega_0^{n-1} )= 0$ are Ricci-flat K\"ahler metrics.
\end{corollary}

\noindent {\it Proof:} The case $n=3$ was done in \cite{PPZ5}, so we assume $n \geq 4$. Setting (\ref{metric-evol-n>3}) to zero and taking the trace gives 
\be
0 =  -\tilde{R} + {n \over 2 (n-2)} (|T|^2- 2 |\tau|^2) - {1 \over 2} |T|^2 + 3 |\tau|^2 .
\ee
By the definition of Ricci curvatures (\ref{def-Ric}),
\be
\tilde{R} = g^{j \bar{k}} R^p{}_{p \bar{k} j} = - g^{j \bar{k}} \p_{\bar{k}} \p_j \log \det g_{\bar{\ell} m},
\ee
hence $\tilde{R} = \Delta \log \| \Omega \|^2_\omega$. Simplifying, we obtain
\be
(n-2) \Delta \log \| \Omega \|^2_\omega = |T|^2 + 2(n-3) |\tau|^2.
\ee
By the maximum principle, we conclude that $\log \| \Omega \|^2_\omega$ is constant and $|T|^2=0$. It follows that $\omega$ is K\"ahler with zero Ricci curvature. Q.E.D.

\section{A complex Monge-Amp\`ere flow}

\subsection{The Anomaly flow and a scalar Monge-Amp\`ere flow}
\setcounter{equation}{0}

Let $X$ be a compact K\"ahler manifold with K\"ahler metric $\hat{\chi} = i \hat\chi_{\bar{k} j} dz^j \wedge d \bar{z}^k$. Given a potential function $\varphi: X \rightarrow {\bf R}$, we will use the notation
\be
F[\varphi] = {(\hat{\chi} + i \ddb \varphi)^n \over \hat{\chi}^n}.
\ee
Let $f \in C^\infty(X,{\bf R})$ be a given function. Consider the Monge-Amp\`ere flow of potentials
\be \label{MA-flow2}
\p_t \varphi = e^{-f} F[\varphi], \ \ \varphi(x,0)=0,
\ee
subject to the plurisubharmonicity condition
\be
\chi = \hat{\chi} + i \ddb \varphi >0.
\ee

Assume further that $X$ admits a nowhere holomorphic $(n,0)$-form $\Omega$ and consider the flow (\ref{AF}). We claim that, with the function $f$ chosen by
\be
e^{-f} = {1 \over (n-1)}  \| \Omega \|^{-2}_{\hat{\chi}},
\ee
the Monge-Amp\`ere flow (\ref{MA-flow2}) is just the Anomaly flow (\ref{AF}) with initial data
\be\label{initialansatz}
\|\Omega\|_{\o_0}\o_0^{n-1}= \hat{\chi}^{n-1}.
\ee 
Indeed, let $t\to \varphi(t)$ be the solution of the Monge-Amp\`ere flow (\ref{MA-flow2}) and set
\be \label{ansatz}
\chi(t) = \hat{\chi} + i \ddb \varphi(t) > 0, \ \ \| \Omega \|_{\omega(t)} \omega(t)^{n-1} = \chi(t)^{n-1}.
\ee
Then 
\be
\p_t ( \| \Omega \|_\omega \omega^{n-1}) = \p_t \chi^{n-1} = (n-1)  i \ddb (\p_t \varphi) \wedge \chi^{n-2}.
\ee
Next, the above relation between $\o(t)$ and $\chi(t)$ can be inverted: taking the determinants of both sides gives $\| \Omega \|_\omega^{1/(n-1)} = \| \Omega \|_\chi^{2/(n-2)}$, and thus the relation can also be written as
\be  \label{omega-ansatz-solved}
\omega^{n-2} = \| \Omega \|_\chi^{-2} \chi^{n-2}.
\ee 
It follows that
\be
i \ddb \omega^{n-2} = i \ddb ( \| \Omega \|_\chi^{-2} \chi^{n-2} ) = i \ddb (\| \Omega \|_\chi^{-2}) \wedge \chi^{n-2}.
\ee
It suffices then to show that
\be
(n-1)\p_t\varphi=\|\Omega\|_\chi^{-2}
\ee
since this would imply that $\o$ is a solution of the Anomaly flow (\ref{AF}), which is uniquely determined by the parabolicity of the flow established in Theorem \ref{main1}. But by the definition of $\varphi$, we have
\bea
(n-1)\p_t\varphi=(n-1)
e^{-f}F[\varphi]={{\rm det}\,\hat\chi\over \Omega\bar\Omega}\cdot
{{\rm det}\,(\hat\chi+i\ddb\varphi)\over{\rm det}\,\hat\chi}
=
{1\over \|\Omega\|_\chi^2}.
\eea
This completes the proof of our claim.

\subsection{Evolution identities}
In this section, we compute the evolution of various basic quantities along the flow (\ref{MA-flow}). Let us first establish notations.
\smallskip
\par We will denote $\chi_{\bar{k} j} = \hat{\chi}_{\bar{k} j} + \varphi_{\bar{k} j}$, and use upper indices $\chi^{j \bar{k}}$ to denote the inverse matrix of $\chi_{\bar{k} j}$. We will sometimes write $\dot{\varphi}$ for $\p_t \varphi$. Subscripts on a function will denote spacial partial derivatives and $\hat{\na}$ and $\Delta_{\hat{\chi}}=\hat{\chi}^{p \bar{q}} \hat{\na}_p \hat{\na}_{\bar{q}}$ will denote covariant derivatives with respect to the background K\"ahler metric $\hat{\chi}$. We use the Chern connection defined as usual by $\hat{\na}_{\bar{k}} V^p = \p_{\bar{k}} V^p$, $\hat{\na}_k V^p = \hat{\chi}^{p \bar{q}} \p_k (\hat{\chi}_{\bar{q} i} V^i)$ for any section $V$ of $T^{1,0}X$, and the curvature convention
\be
[\hat{\na}_{j}, \hat{\na}_{\bar{k}} ] W_i = -\hat{R}_{\bar{k} j}{}^p{}_{i} W_p,
\ee
for any section $W$ of $(T^{1,0}X)^*$. We introduce the linearized operator
\bea\label{defL}
L = e^{-f} F \, \chi^{j \bar{k}} \hat{\na}_j \hat{\na}_{\bar{k}}.
\eea
We will also use the relative endomorphism $h^i{}_j = \hat{\chi}^{i \bar{k}} \chi_{\bar{k} j}$, which satisfies
\be
{\rm Tr} \, h = n + \Delta_{\hat{\chi}} \varphi, \ \ {\rm Tr} \, h^{-1} = \chi^{p \bar{q}} \hat{\chi}_{\bar{q} p}.
\ee

\subsubsection{Evolution of the potential}
We first compute
\bea
(\p_t - L) \varphi &=&  e^{-f} F  - e^{-f} F \chi^{k \bar{j}} \varphi_{\bar{j} k} \nonumber\\
&=&  e^{-f} F - e^{-f} F \chi^{k \bar{j}} (\chi_{\bar{j} k} - \hat{\chi}_{\bar{j} k}).
\eea
Therefore
\be \label{evol_phi}
(\p_t - L) \varphi = - (n-1) e^{-f}F  +  e^{-f} F \, {\rm Tr}\, h^{-1}.
\ee

\subsubsection{Evolution of the determinant}
Next, differentiating the evolution equation in time gives
\be \label{evol_F_0}
\p_t (e^{-f}F) = L (e^{-f} F).
\ee
Expanding, we obtain
\be \label{evol_F}
(\p_t - L) F = e^f L(e^{-f}) F -2  F e^{-f} \, \Re \{ \chi^{j \bar{k}} f_j F_{\bar{k}} \}.
\ee

\subsubsection{Evolution of the cohomological representative} 
Differentiating equation (\ref{MA-flow}) once gives
\be \label{MA-flow-diff}
\p_t \p_{\bar{p}} \varphi = e^{-f} F \chi^{k \bar{j}} \hat{\na}_{\bar{p}} \varphi_{\bar{j} k}  - e^{-f} F \, \p_{\bar{p}} f.
\ee
We introduce the notation
\be \label{F^kjrs}
F^{k \bar{j}, r \bar{s}} = F \chi^{r \bar{s}} \chi^{k \bar{j}} - F \chi^{r \bar{j}} \chi^{k \bar{s}},
\ee
so that
\be
\hat{\na}_q (F \chi^{k \bar{j}} ) = F^{k \bar{j}, r \bar{s}} \hat{\na}_q \varphi_{\bar{s} r}.
\ee
We note that unlike when $F$ is concave, the quantity $-F^{k \bar{j}, r \bar{s}} M_{\bar{j} k} \overline{M_{\bar{r} s}}$ for a Hermitian tensor $M_{\bar{a} b}$ does not have a sign. Differentiating the flow twice, which amounts to differentiating (\ref{MA-flow-diff}) once more gives
\bea
\p_t \varphi_{\bar{p} q} &=& e^{-f} F \chi^{k \bar{j}} \hat{\na}_q \hat{\na}_{\bar{p}} \varphi_{\bar{j} k} + e^{-f} F^{k \bar{j}, r \bar{s}} \hat{\na}_{\bar{p}} \varphi_{\bar{j} k} \, \hat{\na}_{q} \varphi_{\bar{s} r}   \nonumber\\
&& - e^{-f} F \chi^{k \bar{j}} \hat{\na}_{\bar{p}} \varphi_{\bar{j} k} \p_q f  - e^{-f} F \, \p_q \p_{\bar{p}} f  + e^{-f} F \, \p_{\bar{p}} f \p_q f\nonumber\\
&& - e^{-f} F \chi^{k \bar{j}} \hat{\na}_q \varphi_{\bar{j} k} \, \p_{\bar{p}} f.
\eea
Commuting derivatives
\be
e^{-f} F \chi^{k \bar{j}} \hat{\na}_q \hat{\na}_{\bar{p}} \varphi_{\bar{j} k}  = e^{-f} F \chi^{k \bar{j}} \hat{\na}_k \hat{\na}_{\bar{j}} \varphi_{\bar{p} q} - e^{-f}F \chi^{k \bar{j}} \hat{R}_{\bar{j} k \bar{p}}{}^{\bar{a}} \varphi_{\bar{a} q} + e^{-f} F \chi^{k \bar{j}} \hat{R}_{\bar{j} q \bar{p}}{}^{\bar{a}} \varphi_{\bar{a} k}.
\ee
Thus we obtain the evolution of $\chi_{\bar{p} q} = \hat{\chi}_{\bar{p} q} + \varphi_{\bar{p} q}$
\bea \label{evol_chi}
(\p_t -L)\chi_{\bar{p} q} &=& -e^{-f} F \chi^{k \bar{j}} \hat{R}_{\bar{j} k \bar{p}}{}^{\bar{a}} \varphi_{\bar{a} q} + e^{-f} F \chi^{k \bar{j}} \hat{R}_{\bar{j} q \bar{p}}{}^{\bar{a}} \varphi_{\bar{a} k} + e^{-f} F^{k \bar{j}, r \bar{s}} \hat{\na}_{\bar{p}} \chi_{\bar{j} k} \hat{\na}_{q} \chi_{\bar{s} r}   \nonumber\\
&& - e^{-f} F \chi^{k \bar{j}} \hat{\na}_{\bar{p}} \chi_{\bar{j} k} \, f_q  -  e^{-f} F \chi^{k \bar{j}} \hat{\na}_{q} \chi_{\bar{j} k} \, f_{\bar{p}} \nonumber\\
&&- e^{-f} F \, f_{\bar{p} q}   + e^{-f} F \, f_{\bar{p}}  f_q. 
\eea
Since $L$ involves covariant derivatives with respect to the background $\hat{\chi}$, we may take the trace with respect to $\hat{\chi}$ and derive
\bea \label{evol_trh}
(\p_t -L) {\rm Tr}\, h &=& -e^{-f} F \chi^{k \bar{j}} \hat{R}_{\bar{j} k}{}^{p\bar{a}} \varphi_{\bar{a} p} + e^{-f} F \chi^{k \bar{j}}  \hat{R}_{\bar{j}p}{}^{p\bar{a}} \varphi_{\bar{a} k}   \nonumber\\
&& + e^{-f} F^{k \bar{j}, r \bar{s}}  \hat{\chi}^{p \bar{q}} \hat{\na}_{\bar{q}}  \chi_{\bar{j} k} \hat{\na}_{p} \chi_{\bar{s} r} - 2 e^{-f} \Re \{ \hat{\chi}^{q \bar{p}} F_q f_{\bar{p}} \} \nonumber\\
&& - e^{-f} F \hat{\chi}^{q \bar{p}} f_{\bar{p} q}   + e^{-f} F \hat{\chi}^{q \bar{p}} f_q f_{\bar{p}}  . 
\eea
Writing $\chi_{\bar{k} \ell} = \hat{\chi}_{\bar{k} \ell} + \varphi_{\bar{k} \ell}$ and $\hat{R}^p{}_p{}^q{}_q=\hat{R}$,
\bea \label{evol_trh}
(\p_t -L) {\rm Tr}\, h &=& -e^{-f} F \chi^{k \bar{j}} \hat{R}_{\bar{j} k}{}^{p\bar{a}} \chi_{\bar{a} p} + e^{-f} F \hat{R}  + e^{-f} F^{k \bar{j}, r \bar{s}}  \hat{\chi}^{p \bar{q}} \hat{\na}_{\bar{q}}  \chi_{\bar{j} k} \hat{\na}_{p} \chi_{\bar{s} r} \nonumber\\
&&- 2 e^{-f} \Re \{ \hat{\chi}^{q \bar{p}} F_q f_{\bar{p}} \}  - e^{-f} F \hat{\chi}^{q \bar{p}} f_{\bar{p} q}   + e^{-f} F \hat{\chi}^{q \bar{p}} f_q f_{\bar{p}}  . 
\eea

\subsection{Estimate of the speed}
Differentiating the equation in time, we obtain
\be
(\p_t - L) \dot{\varphi} = 0.
\ee
By the maximum principle, we have the bound
\be
\inf_X \dot{\varphi} (0) \leq \dot{\varphi} \leq \sup_X \dot{\varphi} (0).
\ee
Since $\dot{\varphi} (0) = e^{-f}$, it follows that $e^{-f} F = \dot{\varphi} \geq \inf_X e^{-f}$. As a consequence, we obtain
\begin{lemma} \label{C0-est}
Let $\varphi$ be a solution to (\ref{MA-flow}) on $X \times [0,T)$ satisfying the positivity condition (\ref{positivity}). There exists a constants $C>0$ and $\delta>0$ depending only on $(X,\hat{\chi})$ and $f$ such that
\be
|\p_t \varphi | \leq C, \ \ \delta \leq e^{-f} F[\varphi] \leq C, \ \ (\sup_{X} \varphi - \inf_{X} \varphi)(t) \leq C.
\ee
\end{lemma}

{\it Proof:} It only remains to establish the oscillation estimate. This follows from applying Yau's $C^0$ estimate \cite{Y} to the equation
\be
(\hat{\chi} + i \ddb \varphi)^n = ( \p_t \varphi ) \, e^f \, (\hat{\chi})^n
\ee
at each fixed time. Q.E.D.

\subsection{Second order estimate}
\begin{lemma}\label{C2-est}
Let $\varphi$ be a solution to (\ref{MA-flow}) on $X \times [0,T]$ satisfying the positivity condition (\ref{positivity}). There exists a constants $C>0$ and $A>0$ depending only on $X$, $\hat{\chi}$ and $f$ such that
\be
\Delta_{\hat{\chi}} \, \varphi(x,t) \leq C e^{A (\tilde{\varphi}(x,t) - \inf_{X \times [0,T]} \tilde{\varphi})},
\ee
where
\be
\tilde{\varphi}(x,t) = \varphi - {1 \over V} \int_X \varphi \, \hat{\chi}^n, \ \ V = \int_X \hat{\chi}^n.
\ee
\end{lemma}
We will apply the maximum principle to the following test function defined on $X \times [0,T]$,
\be
G(x,t) = \log {\rm Tr}\, h - A \tilde{\varphi} + {B \over 2} F^2,
\ee
where $A,B>0$ are constants to be determined. Before preceeding, we note the extra term ${B\over 2}F^2$, which does not appear in the standard test function used in the study of the complex Monge-Amp\`ere type equations (see for example, \cite{PPZ3, Y}) or the K\"ahler-Ricci flow \cite{Cao}. Indeed, this is the main innovation of this test function, and we use it to overcome the difficult terms caused by the lack of concavity of the Monge-Amp\`ere flow, as well as the cross terms involving the conformal factor $e^{-f}$.

\medskip

We compute the evolution of the test function.
\bea
(\p_t - L) G &=& {1 \over {\rm Tr}\, h} (\p_t - L) {\rm Tr} \, h + {e^{-f} F \over ({\rm Tr}\, h)^2} \chi^{j \bar{k}} (\p_j {\rm Tr}\, h)( \p_{\bar{k}} {\rm Tr}\, h)  - A (\p_t - L) \varphi \nonumber\\
&& + {A \over V} \int_X \p_t \varphi \, \hat{\chi}^n + B F (\p_t - L) F - B e^{-f}F \chi^{j \bar{k}} F_j F_{\bar{k}}.
\eea
Substituting (\ref{evol_phi}), (\ref{evol_F}) and (\ref{evol_trh})
\bea 
(\p_t - L) G &=& {1 \over {\rm Tr}\, h} \bigg\{ -e^{-f} F \chi^{k \bar{j}} \hat{R}_{\bar{j} k}{}^{p \bar{a}} \chi_{\bar{a} p} + e^{-f} F \hat{R} - e^{-f} F \hat{\chi}^{j \bar{k}} f_{\bar{k} j}   + e^{-f} F \hat{\chi}^{j \bar{k}} f_j f_{\bar{k}}   \nonumber\\
&&+ e^{-f} F^{k \bar{j}, r \bar{s}}  \hat{\chi}^{p \bar{q}} \hat{\na}_{\bar{q}}  \chi_{\bar{j} k} \hat{\na}_{p} \chi_{\bar{s} r} - 2 e^{-f}  \Re \{ \hat{\chi}^{j \bar{k}} F_j  f_{\bar{k}} \}   + {e^{-f} F \over {\rm Tr}\, h} \chi^{j \bar{k}} (\p_j {\rm Tr}\, h)( \p_{\bar{k}} {\rm Tr}\, h)\bigg\}  \nonumber\\
&& + A (n-1) e^{-f}F  -A e^{-f} F \, {\rm Tr}\, h^{-1} + {A \over V} \int_X e^{-f} F \, \hat{\chi}^n \nonumber\\
&&+ B F^3 \chi^{j \bar{k}} (e^{-f})_{\bar{k} j} -2 B F^2 e^{-f} \, \Re \{ \chi^{j \bar{k}} f_j F_{\bar{k}} \} - B e^{-f} F \chi^{j \bar{k}} F_j F_{\bar{k}}.
\eea
Using Lemma \ref{C0-est}, we estimate
\bea
 &\ & {1 \over {\rm Tr}\, h} \bigg\{ -e^{-f} F \chi^{k \bar{j}} \hat{R}_{\bar{j} k}{}^{p \bar{a}} \chi_{\bar{a} p} + e^{-f} F \hat{R}  - e^{-f} F \hat{\chi}^{j \bar{k}} f_{\bar{k} j}  + e^{-f} F \hat{\chi}^{j \bar{k}} f_j f_{\bar{k}} \bigg\} \nonumber\\
&\leq& {C \over {\rm Tr}\, h} \{ ({\rm Tr} \, h^{-1})({\rm Tr} \, h)  + 1 \}\leq C {\rm Tr} \, h^{-1},
\eea
where the constant $C$ depends on $e^{-f}$ and the curvature of the background metric $\hat{\chi}$. Next, by (\ref{F^kjrs}) we have
\bea
 e^{-f} F^{k \bar{j}, r \bar{s}}  \hat{\chi}^{p \bar{q}} \hat{\na}_{\bar{q}}  \chi_{\bar{j} k} \hat{\na}_{p} \chi_{\bar{s} r} &=&  e^{-f} F \chi^{r \bar{s}} \chi^{k \bar{j}}  \hat{\chi}^{p \bar{q}} \hat{\na}_{\bar{q}}  \chi_{\bar{j} k} \hat{\na}_{p} \chi_{\bar{s} r} - e^{-f} F \chi^{r \bar{j}} \chi^{k \bar{s}} \hat{\chi}^{p \bar{q}} \hat{\na}_{\bar{q}}  \chi_{\bar{j} k} \hat{\na}_{p} \chi_{\bar{s} r} \nonumber\\
&=& {e^{-f} \over F} \hat{\chi}^{p \bar{q}} F_p F_{\bar{q}} - e^{-f} F  \chi^{r \bar{j}} \chi^{k \bar{s}} \hat{\chi}^{p \bar{q}} \hat{\na}_{\bar{q}}  \chi_{\bar{j} k} \hat{\na}_{p} \chi_{\bar{s} r}.
\eea
The term involving $\hat{\chi}^{p \bar{q}} F_p F_{\bar{q}}$ is the new bad term compared to standard arguments, and it is the reason for the addition of ${B \over 2} F^2$ to the test function $G$. The main inequality becomes
\bea 
(\p_t - L) G &\leq& {1 \over {\rm Tr}\, h} \bigg\{ {e^{-f} \over F} \hat{\chi}^{j \bar{k}} F_j F_{\bar{k}}  - 2 e^{-f}  \Re \{ \hat{\chi}^{j \bar{k}} F_j f_{\bar{k}} \} \bigg\} \nonumber\\
&& +  {e^{-f} F \over {\rm Tr} \, h} \bigg\{{1 \over {\rm Tr}\, h} \chi^{j \bar{k}} (\p_j {\rm Tr}\, h)( \p_{\bar{k}} {\rm Tr}\, h)  -   \hat{\chi}^{p \bar{q}} \chi^{r \bar{j}} \chi^{k \bar{s}} \hat{\na}_{\bar{q}}  \chi_{\bar{j} k} \hat{\na}_{p} \chi_{\bar{s} r}  \bigg\} \nonumber\\
&& -2 B F^2 e^{-f} \, \Re \{ \chi^{j \bar{k}}f_j F_{\bar{k}} \} - B e^{-f} F \chi^{j \bar{k}} F_j F_{\bar{k}} \nonumber\\
&&   + \{ C(1+B)-A e^{-f} F \} {\rm Tr} \, h^{-1} + C(A+B+1).
\eea
By the inequality of Yau and Aubin \cite{Y,Au}, 
\be
{1 \over {\rm Tr}\, h} \chi^{j \bar{k}} (\p_j {\rm Tr}\, h)( \p_{\bar{k}} {\rm Tr}\, h) -  \hat{\chi}^{p \bar{q}} \chi^{r \bar{j}} \chi^{k \bar{s}}  \hat{\na}_{\bar{q}}  \chi_{\bar{j} k} \hat{\na}_{p} \chi_{\bar{s} r}  \leq 0.
\ee
Next, we estimate
\bea
{1 \over {\rm Tr}\, h} \bigg\{ {e^{-f} \over F} \hat{\chi}^{j \bar{k}} F_j F_{\bar{k}}  - 2 e^{-f} \Re \{ \hat{\chi}^{j \bar{k}} F_j f_{\bar{k}} \} \bigg\} &\leq& {1 \over {\rm Tr}\, h} \bigg\{ 2 {e^{-f} \over F}\hat{\chi}^{j \bar{k}} F_j F_{\bar{k}} + e^{-f} F \hat{\chi}^{j \bar{k}} f_j f_{\bar{k}} \bigg\} \nonumber\\
&\leq&  2 {e^{-f} \over F} \chi^{j \bar{k}} F_j F_{\bar{k}} + C {\rm Tr}\, h^{-1},
\eea
and
\be
-2 B F^2 e^{-f} \, \Re \{ \chi^{j \bar{k}}f_j F_{\bar{k}} \} \leq  \chi^{j \bar{k}} F_j F_{\bar{k}} + CB^2 {\rm Tr} \, h^{-1}.
\ee
Applying these estimates in the main inequality yields
\bea 
(\p_t - L) G &\leq& (1+2 e^{-f} F^{-1} - B  e^{-f} F) \chi^{j \bar{k}} F_j F_{\bar{k}} \nonumber\\
&& + \{ C(1+B+B^2)-A e^{-f} F\} {\rm Tr} \, h^{-1} + C(A+B+1).
\eea
By Lemma \ref{C0-est}, we have the uniform lower bound $e^{-f} F \geq \delta>0$. Thus may choose $B \gg 1$ such that $(1+2 e^{-f} F^{-1} - B  e^{-f} F)\leq 0$. Next, we may choose $A \gg B \gg 1$ such that $(C(1+B+B^2)-A e^{-f} F) \leq -1$. Then if $G$ attains a maximum on $X \times [0,T]$ at a point $(x_0,t_0)$ with $t_0>0$, by the maximum principle we have the inequality
\be
0 \leq (\p_t - L)G \leq - {\rm Tr} \, h^{-1} + C(A+B+1).
\ee
It follows that the eigenvalues of $h$ are bounded below at $(x_0,t_0)$. The product of the eigenvalues of $h$ is given by $F$, which is uniformly bounded along the flow. Thus the eigenvalues of $h$ are bounded above at $(x,t_0)$, and so ${\rm Tr} \, h \leq C$ at $(x_0,t_0)$. Therefore
\be \label{C2-G-est}
G(x,t) \leq G(x_0,t_0) \leq C - A \inf_{X \times [0,T]} \tilde{\varphi}.
\ee
If $G$ attains a maximum at $t_0=0$, we have already have $G(x,t) \leq C$. Therefore 
\be
\log {\rm Tr} \, h \leq C+A \left( \tilde{\varphi} - \inf_{X \times [0,T]} \tilde{\varphi} \right)
\ee
 along the flow. Q.E.D.

\begin{corollary} \label{unif-ellipt}
Let $\varphi$ be a solution to (\ref{MA-flow}) on $X \times [0,T]$ satisfying the positivity condition (\ref{positivity}). There exists a constant $C>0$ depending only on $X$, $\hat{\chi}$ and $f$ such that
\be
C^{-1} \hat{\chi}_{\bar{k} j}(z) \leq \chi_{\bar{k} j}(z,t) \leq C \hat{\chi}_{\bar{k} j}(z).
\ee
\end{corollary}

\noindent {\it Proof:} We first use Lemma \ref{C2-est} to obtain an upper bound on $\Delta_{\hat{\chi}} \varphi$ independent of time. Indeed, the oscillation of $(\sup_X \tilde{\varphi} - \inf_X \tilde{\varphi})(t)$ is bounded at each fixed time by Lemma \ref{C0-est}. Since $\int \tilde{\varphi} \, \hat{\chi}^n = 0$, we must have that $(\sup_X \tilde{\varphi})(t) \geq 0$ and $(\inf_X \tilde{\varphi})(t) \leq 0$, hence $\tilde{\varphi}$ is uniformly bounded along the flow. By Lemma \ref{C2-est},
\be
{\rm Tr} \, h \leq C.
\ee
This gives the upper bound of $\chi$. For the lower bound, we use that the determinant $F[\varphi] = \chi^n / \hat{\chi}^n$ is uniformly bounded below by Lemma \ref{C0-est}. Q.E.D.

\subsection{Higher order estimates}

In this section, we will follow the main line of the argument given in Tsai \cite{Tsai} to prove the higher order estimates. We note that we cannot directly apply a standard theorem to our equation $\p_t \varphi = e^{-f} F[\varphi]$, $F[\varphi] = \det (\hat{\chi}_{\bar{k} j} + \varphi_{\bar{k} j}) / \det \hat{\chi}_{\bar{k} j}$, since $F[\varphi]$ is not concave in the second derivatives of $\varphi$. However, $\log F[\varphi]$ is concave as a function of $\varphi_{\bar{k} j}$, and to introduce the logarithm we split the argument into several steps to treat space and time seperately. First, we fix some notation.  

\par We will work locally on a ball $B_r(0)$ and cylinder $Q= B_r \times (T_0,T)$. For functions $v: B_r \rightarrow {\bf R}$ and $u: Q \rightarrow {\bf R}$, we define
\be
\| v \|_{C^\alpha(B_r)} = \| v \|_{L^\infty(B_r)} + \sup_{x \neq y \in B_r} {|v(x)-v(y)| \over |x-y|^\alpha},
\ee
and
\be
\| u \|_{C^{\alpha,\alpha/2}(Q)} = \| u \|_{L^\infty(Q)} + \sup_{(x,t) \neq (y,s) \in Q} {|u(x,t)-u(y,s)| \over (|x-y|+|t-s|^{1/2})^\alpha}.
\ee
The main estimate of this section is the following.

\begin{lemma} \label{c2-alpha-estimate}
Let $\varphi$ be a solution to (\ref{MA-flow}) on $X \times [0,\epsilon)$ satisfying the positivity condition (\ref{positivity}). Let $B_1$ be a coordinate chart on $X$ such that $B_1 \subset {\bf R}^n$ is a unit ball. Then there exists $0<\alpha<1$ and $C>0$, depending on $\hat{\chi}$, $f$, and $\epsilon$, such that on $Q=B_{1/2} \times [{\epsilon \over 2},\epsilon)$,
 \bea
\| \p_t \varphi \|_{C^{\alpha,\alpha/2}(Q)} + \| \varphi_{\bar{k} j} \|_{C^{\alpha,\alpha/2}(Q)} \leq C.
 \eea
\end{lemma}

We first notice that the speed function of the flow satisfies the equation
\be
\p_t (e^{-f}F) = L(e^{-f}F)=e^{-f}F \chi^{j\bar k} (e^{-f} F)_{\bar k j}.
\ee
By Lemma \ref{C0-est} and Corollary \ref{unif-ellipt}, we see that this is a uniformly parabolic linear equation with bounded coefficients on $X\times [0, \epsilon)$. It follows from the Krylov-Safanov Harnack inequality \cite{KrSa} that
\be
\| e^{-f} F \|_{C^{\alpha,\alpha/2}(Q)} \leq C,
\ee
for some $\alpha\in (0, 1)$. This implies that
\be \label{Holder-speed}
\| F \|_{C^{\alpha,\alpha/2}(Q)} \leq C, \ \ \ \| \p_t \varphi \|_{C^{\alpha,\alpha/2}(Q)}\leq C.
\ee
For each fixed $t \in [{\epsilon \over 2}, \epsilon)$, by covering the compact manifold $X$ we can control the H\"older norm of $F(\cdot,t)$ on all of $X$, and we have the space norm estimate
\bea\label{space-est1}
\| F(\cdot, t) \|_{C^\alpha(B_1)} \leq C,
\eea
where $C$ is a constant independent of $t\in [{\epsilon \over 2}, \epsilon)$. Recall that $F[\varphi] = \det (\hat{\chi}_{\bar{k} j} + \varphi_{\bar{k} j}) / \det \hat{\chi}_{\bar{k} j}$, so after possibly taking a smaller $0<\alpha<1$, we may apply the estimates in Tosatti-Weinkove-Wang-Yang \cite{TWWY} to establish
\bea\label{space-est2}
\| \varphi_{\bar{k} j}(\cdot, t) \|_{C^\alpha(B_1)} \leq C
\eea
where $C$ is independent of $t\in [{\epsilon \over 2}, \epsilon)$. Therefore, in view of a standard lemma \cite{Kr} in parabolic H\"older spaces which allows us to treat time and space separately, it remains to show that
\bea
\sup_{s, t\in [{\epsilon \over 2}, \epsilon], s\neq t} {|\varphi_{\bar k j} (z, s) - \varphi_{\bar k j}(z, t) | \over |s-t|^{\alpha\over 2} } \leq C,
\eea
for some $C$ independent of $z \in B_{1/2}$. Following Tsai \cite{Tsai}, for $0<h<{\epsilon \over 2}$ we consider the function
\bea
U_{\lambda, h}(z, t) = \lambda \varphi(z, t) + (1-\lambda) \varphi(z, t+h)  \textit{ with }0 \leq \lambda \leq 1,
\eea
defined on $B_{1} \times [{\epsilon \over 2}, \epsilon-h)$. Compute
\bea
\log F(z, t) - \log F (z, t+h)
&=&
 \int^1_0 {d\over d\lambda} \log\det\left(\hat\chi_{\bar k j} + \p_{\bar k} \p_j U_{\lambda, h}(z, t) \right) \, d\lambda\\ \nonumber
 &=&
 \left(\int^1_0 \chi_{\lambda, h}^{j\bar k} (z, t) \, d\lambda \right) \cdot ( \varphi_{\bar k j}(z, t) - \varphi_{\bar k j} (z, t+h)) 
\eea
where $\chi_{\lambda, h}^{j\bar k} (z, t)$ is the inverse of $(\hat\chi_{\bar{k} j} + \p_j \p_{\bar{k}} U_{\lambda, h})>0$. Denote 
\be
a_h^{j\bar k} (z, t) =\int^1_0 \chi_{\lambda, h}^{j\bar k} (z, t) \, d\lambda.
\ee
It follows from Corollary \ref{unif-ellipt} and the space norm estimate (\ref{space-est2}) that $a^{j\bar k}$ is uniform elliptic and satisfies
\be
\| a_h^{j\bar k}(\cdot,t) \|_{C^\alpha(B_1)} \leq C
\ee
for some constant $C$ independent of $h$ and $t\in [{\epsilon \over 2}, \epsilon-h)$. Denote
\be
\varphi_h(z, t) = {\varphi(z, t)- \varphi(z, t+h) \over |h|^{\alpha\over 4}}
\ee
 with $(z, t) \in B_1 \times [{\epsilon \over 2}, \epsilon-h)$. Then, $\varphi_h$ satisfies the equation
 \bea\label{equation-uh}
 a_h^{j\bar k} (z, t) \p_j \p_{\bar{k}} \, \varphi_h(z, t) = g_h(z, t)
 \eea
 with 
\be
 g_h(z, t)={\log F(z, t) - \log F (z, t+h)\over |h|^{\alpha\over 4}}.
 \ee
 As we discussed above, at a fixed time we have that $a^{j\bar k}_h( \cdot , t)$ is uniformly elliptic and H\"older continuous, with constants independent of time $t$ and parameter $h$. We need the following lemma to estimate the H\"older norm of $g_h$.

\begin{lemma} \label{g_h-lemma}
The function $g_h(z,t)$ satisfies
\be
\| g_h (\cdot,t) \|_{C^{\alpha \over 4}(B_1)} \leq C,
\ee
where $C$ is independent of $h$ and $t \in [{\epsilon \over 2}, \epsilon - h)$. 
\end{lemma}
Given this lemma, we can apply the elliptic Schauder estimate to equation (\ref{equation-uh}) at a fixed time $t$ and obtain
\be
\| \varphi_h (\cdot, t) \|_{C^2(B_{1/2})} \leq C \left( \| g_h (\cdot,t) \|_{C^{\alpha \over 4}(B_1)} + \| \varphi_h (\cdot, t) \|_{L^\infty(B_1)} \right),
\ee
with $C$ independent of $t$ and $h$. By the estimate of the speed function in Lemma \ref{C0-est}, we have $\| \varphi_h (\cdot, t) \|_{L^\infty(B_1)}\leq C$ with $C$ independent of $t$ and $h$. This together with Lemma \ref{g_h-lemma} implies
\be
\sup_{z\in B_{1/2}} {| \varphi_{\bar{k} j}(z, t) - \varphi_{\bar{k} j}(z, t+h)| \over |h|^{\alpha\over 4} }  \leq C
\ee
for all $h>0$ small, where $C$ is independent of $h$ and $t\in [{\epsilon \over 2}, \epsilon-h)$. Hence, we complete the proof of Lemma \ref{c2-alpha-estimate}.
\bigskip
\par {\it Proof of Lemma \ref{g_h-lemma}:} Since $F[\varphi]$ is uniformly bounded away from zero along the flow, we know that by (\ref{Holder-speed}) that
\be
\| \log F \|_{C^{\alpha,\alpha/2}(B_1 \times [{\epsilon \over 2}, \epsilon))} \leq C.
\ee
This implies that
\be \label{logF-h-ineq}
{1 \over |h|^{\alpha \over 4}} | \log F(x,t) - \log F(x,t+h) - \log F(y,t) + \log F(y,t+h)|
\leq C |x-y|^{\alpha \over 4}.
\ee
Indeed, if $|h| \leq |x-y|$, then
\bea
 &\ &{1 \over |h|^{\alpha \over 4}} \{ |\log F(x,t) - \log F(x,t+h)| + | \log F(y,t) - \log F(y,t+h)| \}  \nonumber\\
&\leq&
C |h|^{\alpha \over 4} \leq C |x-y|^{\alpha \over 4}.
\eea
On the other hand, if $|h| \geq |x-y|$, then
\bea
& \ & {1 \over |h|^{\alpha \over 4}} \{ |\log F(x,t) -  \log F(y,t)| + |\log F(x,t+h) - \log F(y,t+h)| \} \nonumber\\
&\leq& C {|x-y|^{\alpha \over 4} \over |h|^{\alpha \over 4}} |x-y|^{\alpha \over 4} \leq C |x-y|^{\alpha \over 4}.
\eea
Combining these two cases and using the triangle inequality proves (\ref{logF-h-ineq}). Thus for any $t \in [{\epsilon \over 2}, \epsilon - h)$ there holds
\be
 \Big\| {\log F(\cdot, t) - \log F (\cdot, t+h)\over |h|^{\alpha\over 4}}  \Big\|_{C^{\alpha/4}(B_1)} \leq C,
\ee
where $C$ is independent of $t$ and $h$. This completes the proof of the lemma. Q.E.D.
\

\section{Proof of Theorem \ref{main2} and Theorem \ref{main3}}
\setcounter{equation}{0}
 
 In this section, we prove the long time existence and convergence of the complex Monge-Amp\`ere flow, which completes the proof for Theorem \ref{main3}. Then, from the discussion in \S 4.1, we also obtain the proof of Theorem \ref{main2}.
 
\subsection{Long-time existence}

In this section, we show that the solution $\varphi$ and its normalization $\tilde\varphi$ are smooth and exist for all time.

Since the flow is parabolic, a solution exists on a maximal time interval $[0,T)$ with $T>0$. Differentiating the equation, we obtain
\be
\p_t \p_p \varphi = e^{-f} F \chi^{k \bar{j}} \hat{\na}_k \hat{\na}_{\bar{j}} \p_p \varphi - e^{-f} F \p_p \varphi.
\ee
By Corollary \ref{unif-ellipt} and Lemma \ref{c2-alpha-estimate}, this is a uniformly parabolic equation for $\p_p \varphi$ with H\"older continuous coefficients. By parabolic Schauder estimates (e.g. \cite{Kr}), we obtain the $C^{2+\alpha, 1+ \alpha/2}$ estimate for $\p_p \varphi$. A bootstrapping argument gives estimates on all derivatives of $\varphi$. 

If $T<\infty$, then our estimates allow us to take a subsequential limit and then extend the flow using the short-time existence theorem. It follows that a smooth solution exists on $[0,\infty)$. We already noted in the proof of Corollary \ref{unif-ellipt} that $\tilde{\varphi}$ is uniformly bounded, and the above argument shows that $\tilde{\varphi}$ and all its derivatives are bounded along the flow.

\subsection{Dilaton functional}
We now return to the Anomaly flow $\p_t ( \| \Omega \|_\omega \omega^{n-1} ) = i \ddb (\omega^{n-2})$ with ansatz $\| \Omega \|_{\omega} \omega^{n-1} = \chi^{n-1}$. Recall that this ansatz is preserved and the conformally balanced metric $\omega(t)$ is given by the expression 
\be
\omega = \| \Omega \|_\chi^{-2/(n-2)} \chi, \ \ \ \chi(t) = \hat{\chi} + i \ddb \varphi(t),
\ee
where $\varphi$ solves the Monge-Amp\`ere flow (\ref{MA-flow}). By the previous section, the Anomaly flow for the Hermitian metric $\omega(t)$ exists for all time. Let
\be
M(\omega) = \int_X \| \Omega \|_\omega \, \omega^n
\ee
denote the dilaton functional. This functional was introduced by Garcia-Fernandez, Rubio, Shahbazi and Tipler \cite{GFRST} to formulate a variational principle for the Hull-Strominger system on holomorphic Courant algebroids. Here we compute its evolution along the Anomaly flow.

\begin{lemma} \label{evol-dilaton-func}
Let $\omega(t)$ be a solution to the Anomaly flow (\ref{AF}) with $n \geq 3$. Then the dilaton functional evolves by
\be
{d \over dt} M(\omega(t)) = {1 \over 2} {1 \over (n-1)(n-2)} \int_X \{ |T|^2 - 2 |\tau|^2 \} \, \omega^n .
\ee
\end{lemma}
{\it Proof:} Wedging the equation ${d \over dt} (\| \Omega \|_\omega \omega^{n-1}) = i \ddb (\omega^{n-2})$ with $\omega$ gives
\be
\left( {d \over dt} \| \Omega \|_\omega \right) \, \omega^n + (n-1) \| \Omega \|_\omega \omega^{n-1} \wedge \dot{\omega} = i \ddb \omega^{n-2} \wedge \omega. 
\ee
By definition $\| \Omega \|^2_\omega = \Omega \bar{\Omega} (\det g)^{-1}$, and we have
\be
\left( {d \over dt} \| \Omega \|_\omega \right) \, \omega^n = -{n \over 2} \| \Omega \|_\omega \omega^{n-1} \wedge \dot{\omega}.
\ee
Therefore
\be
{(n-2) \over 2} \| \Omega \|_\omega \omega^{n-1} \wedge \dot{\omega} = i \ddb \omega^{n-2} \wedge \omega.
\ee
It follows that
\bea
{d \over dt} \left( \| \Omega \|_\omega \omega^n \right) &=& {d \over dt} \left( \| \Omega \|_\omega \omega^{n-1} \right) \wedge \omega + \| \Omega \|_\omega \omega^{n-1} \wedge \dot{\omega} = \left( {n \over n-2} \right) i \ddb \omega^{n-2} \wedge \omega.\nonumber
\eea
We now use Stokes theorem to obtain
\be
{d \over dt} M(t) = \int_X {d \over dt} \left( \| \Omega \|_\omega \omega^n \right) = - {n \over n-2} \int_X i \p \omega^{n-2} \wedge \bar{\p} \omega.
\ee
Hence
\be
{d \over dt} M(t) = - n \int_X i \p \omega \wedge \bar{\p} \omega \wedge \omega^{n-3}.
\ee
Applying identity (\ref{T-Tbar-contract}) from the appendix with $T=i \p \omega$ and $\bar{T} = - i \bar{\p} \omega$, we may rewrite the integral as
\be 
- n \int_X i \p \omega \wedge \bar{\p} \omega \wedge \omega^{n-3} = {1 \over 2} {1 \over (n-1)(n-2)} \int_X \{ |T|^2 - 2 |\tau|^2 \} \, \omega^n .
\ee
This completes the proof of the lemma. Q.E.D.

\medskip

\begin{lemma} \label{M-evol-conf-kah}
Let $\omega(t)$ be a solution to the Anomaly flow (\ref{AF}) with $n \geq 3$. Further assume that the solution $\omega(t)$ is conformally K\"ahler. Then the dilaton functional evolves by
\be
{d \over dt} M(\omega(t)) = - {1 \over 2 (n-1)} \int_X |T|^2 \, \omega^n .
\ee
In particular, the dilaton functional is monotone decreasing along the flow.
\end{lemma}
{\it Proof:} Fix a time $t$, and let $\omega(t) = e^\psi \hat{\omega}$ with $d \hat{\omega} = 0$. The components of the torsion of $\omega = i g_{\bar{k} j} dz^j \wedge d \bar{z}^k$ are given by
\be
T_{\bar{k} j \ell} = \p_j \psi \, g_{\bar{k} \ell} - \p_\ell \psi \, g_{\bar{k} j}, \ \ \ T_\ell = - (n-1) \p_\ell \psi.
\ee
Computing the norms of these torsion tensors gives
\be
|T|^2 = 2(n-1) |\na \psi|^2_g, \ \ \ |\tau|^2 = (n-1)^2 |\na \psi|^2_g.
\ee
Therefore $2|\tau|^2 = (n-1) |T|^2$ in this case, and the identity follows from Lemma \ref{evol-dilaton-func}. Q.E.D. 

\medskip

\begin{lemma}
Let $\omega(t)$ be a solution on $[0,\infty)$ to the Anomaly flow (\ref{AF}) with $n \geq 3$ with initial data satisfying
\be
\| \Omega \|_{\omega(0)} \omega(0)^{n-1} = \hat{\chi}^{n-1},
\ee
where $\hat{\chi}$ is a K\"ahler metric. Then $M(\omega(t))$ is monotonically decreasing and furthermore ${d \over dt} M(\omega(t)) \rightarrow 0$ as $t \rightarrow \infty$.
\end{lemma}
{\it Proof:} By the formula for the ansatz (\ref{omega-ansatz-solved}), we have $\omega = \| \Omega \|_\chi^{-2/(n-2)} \chi$, and the metric is conformally K\"ahler $\omega(t) = e^{\psi(t)} \chi(t)$ with
\be
\psi = {1 \over (n-2)} \log \| \Omega \|_\chi^{-2}.
\ee
Furthermore
\be
|T|^2= 2(n-1) |\na \psi|^2_\omega = {2(n-1) \over (n-2)^2} \| \Omega\|_\chi^{2/(n-2)} \left| \na \log \| \Omega \|_\chi^{-2} \right|^2_\chi.
\ee
Also,
\be
\omega^n = \| \Omega \|_\chi^{-2n/(n-2)} \chi^n = \| \Omega \|_\chi^{-2n/(n-2)} \| \Omega \|_\chi^{-2} \| \Omega \|^2_{\hat{\chi}} \, \hat{\chi}^n.
\ee
We may apply Lemma \ref{M-evol-conf-kah} and obtain that $M(t)$ is monotonically decreasing along the flow. Recall that $(n-1)\dot{\varphi} = \| \Omega \|_\chi^{-2}$ gives the flow of the potential. Then
\bea
{d \over dt} M(t) &=& - {1 \over (n-2)^2} \int_X (\| \Omega \|_\chi^{-2})^{2n-3 \over n-2} \left| \na \log \| \Omega \|_\chi^{-2} \right|^2_\chi \, \|\Omega\|_{\hat\chi}^2\,\hat{\chi}^n\\
&=&
- {(n-1)^{2n-3\over n-2} \over (n-2)^2 } \int_X \dot{\varphi}^{1\over n-2} \, |\nabla \dot{\varphi}|_{\chi}^2\, \|\Omega\|_{\hat\chi}^2\,\hat{\chi}^n. \nonumber
\eea
We compute
\bea
{d^2 \over dt^2} M(t) &=& - { (n-1)^{2n-3\over n-2} \over (n-2)^2} \bigg\{ {1\over n-2} \int_X \dot{\varphi}^{-(n-3) \over n-2} \p_t \dot{\varphi}\, | \na \dot{\varphi} |^2_\chi \, \|\Omega\|_{\hat\chi}^2\,\hat{\chi}^n \nonumber\\
&& -\int_X \dot{\varphi}^{1\over n-2} \,  \chi^{j \bar{q}} \chi^{p \bar{k}} \p_j \dot{\varphi}\, \p_{\bar{k}}\dot{\varphi} \,\dot{\varphi}_{\bar{q} p}\, \|\Omega\|_{\hat\chi}^2 \, \hat{\chi}^n \nonumber\\
&&+\int_X \dot{\varphi}^{1\over n-2} \,  \chi^{j\bar k} \left(\p_j\p_t\dot{\varphi} \, \p_{\bar k} \dot{\varphi} + \p_j \dot{\varphi} \, \p_{\bar k} \p_t \dot{\varphi} \right) \, \|\Omega\|_{\hat\chi}^2 \, \hat{\chi}^n \bigg\}.
\eea
By our estimates, $\dot{\varphi}$ and $\chi$ are uniformly bounded above and away from zero, and all space-time derivatives of $\dot{\varphi}$ are bounded. $M(t)$ is uniformly bounded, monotone decreasing, and ${d^2 \over dt^2} M(t)$ is uniformly bounded. It follows that ${d \over dt} M(t) \rightarrow 0$ as $t \rightarrow 0$. Q.E.D.


\subsection{Convergence}
\smallskip
It remains only to show the convergence of the Anomaly flow (\ref{AF}),
which we shall do, using the dilaton functional. Suppose there exists a sequence of times $t_j \rightarrow \infty$ such that $\omega(t_j)$ does not converge to $\omega_\infty$ as given in the theorem. By our estimates and the Arzela-Ascoli theorem, upon taking a subsequence we have that $\omega(t_{j_k})$ converges smoothly to a metric $\omega'_\infty \neq \omega_\infty$. Since ${d \over dt} M(t) \rightarrow 0$, we conclude by Lemma \ref{M-evol-conf-kah} that
\be
\int_X |T(\omega'_\infty)|^2 (\omega'_\infty)^n = 0
\ee 
and hence $\omega'_\infty$ is K\"ahler. By the ansatz (\ref{omega-ansatz-solved}),
\be
\omega'_\infty = \| \Omega \|_{\chi'_\infty}^{-2/(n-2)} \chi'_\infty,
\ee
and so $\| \Omega \|_{\chi'_\infty}$ is constant. It follows that $\chi'_\infty = \hat{\chi} + i \ddb \varphi_\infty$ is the unique K\"ahler Ricci-flat metric in $[\hat{\chi}]$. Since
\be
\int_X \| \Omega \|_{\chi'_\infty}^2 {(\chi'_\infty)^n \over n!} = \int_X i^{n^2} \Omega \wedge \bar{\Omega},
\ee
the constant $\| \Omega \|_{\chi'_\infty}$ is identified. Thus $\omega'_\infty = \omega_\infty$ given in the theorem, and we have smooth convergence. Q.E.D.
\bigskip

\par 
We chose the argument using the dilaton functional in the belief that it will be useful in future studies of the Anomaly flow. For those readers who are only interested in the scalar equation (\ref{MA-flow}), there are alternate ways to establish convergence of $\tilde{\varphi}$. For example, the functional $\int_X e^{-f} F^2 \hat{\chi}^n$ is also monotone decreasing. Alternatively, convergence can be obtained by using the Li-Yau Harnack inequality as in \cite{Cao, Gill, PT}. In this case, we would use the Li-Yau Harnack estimate for the heat equation $u_t = g^{j\bar k} u_{\bar k j}$ on a Hermitian manifold $(M,g)$ proved by Gill \cite{Gill}, and apply it to the differentiated equation $\p_t \dot{\varphi} = e^{-f} F \chi^{j \bar{k}} \dot{\varphi}_{\bar{k} j}$.

\section{Further Developments}
\setcounter{equation}{0}

We conclude with some observations and open questions.

\medskip
(a) The convergence theorems established in this paper for the flow (\ref{AF}) should be viewed as only the first step in a fuller theory yet to be developed. 
For example, we do not know at this moment whether the flow (\ref{AF}) will converge if the initial data is only known to satisfy $\|\Omega\|_{\o_0}\o_0^{n-1}
\in [\hat\chi^{n-1}]$. We expect that it will not, unless $\|\Omega\|_{\o_0}\o_0^{n-1}=(\chi')^{n-1}$ for some K\"ahler form $\chi'$, in which case $[\chi']=[\hat\chi]$. If so, whether the flow (\ref{AF}) converges with initial data $\|\Omega\|_{\o_0}\o_0^{n-1}$ may serve as a criterion for whether $\|\Omega\|_{\o_0}\o_0^{n-1}$ is the $(n-1)$-th power of a K\"ahler form.

In general, we actually expect the failure of convergence of the flow to provide important geometric information. As just stated above, this failure may be caused by the choice of initial data, even within the $(n-1,n-1)$-cohomology class $[\hat\chi^{n-1}]$ with $\hat\chi$ K\"ahler. More important, it may be caused by the class $[\|\Omega\|_{\omega_0}\omega_0^{n-1}]$ not containing $\hat\chi^{n-1}$ for any K\"ahler form $\chi$, and in particular by the manifold $X$ not being K\"ahler. In all these situations, we expect the formation of singularities of the flow and/or long-time behavior to reflect the non-K\"ahler setting. We shall return to these issues elsewhere.

\medskip

(b) The existence of an initial metric $\o_0$ satisfying the condition (\ref{initialansatz}) is equivalent to the existence of a K\"ahler metric $\hat\chi$. Indeed, if $\hat\chi$ is a K\"ahler metric, we can set $\o_0=\|\Omega\|_{\hat\chi}^{-{2\over n-2}}\hat\chi$ to obtain a metric satisfying (\ref{initialansatz}). 

\par In fact, any conformally balanced initial metric which is conformally K\"ahler satisfies (\ref{initialansatz}). Let $\omega_0$ be an initial conformally balanced metric such that $\omega_0 = e^{\psi}\hat{\chi}$ where $\psi: X \rightarrow {\bf R}$ is a smooth function and $\hat{\chi}$ is a K\"ahler metric. Substituting $\omega_0 = e^{\psi} \hat{\chi}$, we obtain
\be
\| \Omega \|_{\omega_0} \omega_0^{n-1} = (e^{({n\over 2}-1)\psi} \| \Omega \|_{\hat{\chi}}) \hat{\chi}^{n-1}.
\ee
Since $d (\| \Omega \|_{\omega_0} \omega_0^{n-1})=0$, we conclude that $e^{({n\over 2}-1)\psi} \| \Omega \|_{\hat{\chi}} = C$ where $C>0$ is a constant. It follows that $\| \Omega \|_{\omega_0} \omega_0^{n-1} = C \hat{\chi}^{n-1}$. After replacing $\hat{\chi}$ by $C^{1/(n-1)} \hat{\chi}$, we see that the ansatz (\ref{initialansatz}) is satisfied. 

\par In particular, we have shown that the Anomaly flow (\ref{AF}) preserves the conformally K\"ahler condition.

\medskip

\par (c) Given an initial metric $\omega_0$ satisfying $d (\| \Omega \|_{\omega_0} \omega_0^{n-1}) = 0$, the Anomaly flow (\ref{AF}) preserves the balanced class of the initial metric.
\be
{d \over dt} [\| \Omega \|_{\omega(t)} \omega(t)^{n-1}] = [i \ddb \omega^{n-2}] = 0.
\ee
Here we take cohomology classes in Bott-Chern cohomology
\be
H_{BC}^{n-1,n-1}(X,{\bf R}) = {\{ {\rm closed} \ {\rm real} \ (n-1,n-1) \ {\rm forms} \} \over \{ i \ddb \beta : \beta \in \Omega^{n-2,n-2}(X,{\bf R}) \} }.
\ee
Thus
\be
\| \Omega \|_{\omega(t)} \omega(t)^{n-1} \in [\| \Omega \|_{\omega_0} \omega_0^{n-1}],
\ee
where
\be
[\| \Omega \|_{\omega_0} \omega_0^{n-1}] \in H_{BC}^{n-1,n-1}(X,{\bf R})
\ee
is the balanced class of $\omega_0$. Since stationary points of the Anomaly flow are K\"ahler metrics, the Anomaly flow could potentially be used to study the relation between the balanced cone and K\"ahler cone on a K\"ahler Calabi-Yau manifold $(X,\Omega)$. The interaction between these two cones was explored by J.-X. Fu and J. Xiao \cite{FX}, and they raised the question of detecting when a balanced class contains a K\"ahler metric (or more generally, a limit of K\"ahler metrics). Examples are given in \cite{FX}, \cite{To} of positive balanced classes on a K\"ahler manifold which do not contain a K\"ahler metric. As discussed above in (a), the formation of singularities of the Anomaly flow may be related to the properties of the initial balanced class. 

\medskip

(d) Similar questions to (a) have been brought to our attention in informal discussions with T. Collins, in connection with a conjecture of Lejmi-Sz\'ekelyhidi \cite{LS}. We briefly describe here one special case of the Lejmi-Sz\'ekelyhidi conjecture. Let $\omega$ and $\alpha$ be two K\"ahler metrics on $X$. It is conjectured that if
\be
\int_X (\omega^n - n \omega \wedge \alpha^{n-1}) \geq 0, \ \  \int_D ( \omega^{n-1} - \alpha^{n-1}) > 0,
\ee
for every irreducible divisor $D$, then there exists a K\"ahler metric $\omega' \in [\omega]$ satisfying the positivity condition
\be \label{subsoln-quotient}
\omega'^{n-1} - \alpha^{n-1} > 0.
\ee
It was proved by Xiao \cite{Xiao2} that one can find a balanced Hermitian metric $\tilde{\omega}$ such that $[\tilde{\omega}^{n-1}] = [\omega^{n-1}]$ and $\tilde{\omega}^{n-1} - \alpha^{n-1}>0$, but it remains to find a K\"ahler metric in the balanced class $[\omega^{n-1}]$ with the desired positivity. The positivity condition (\ref{subsoln-quotient}) is of interest, as it corresponds to subsolutions to the fully nonlinear PDE
\be
n \alpha^{n-1} \wedge (\omega+i \ddb u) = (\omega+i \ddb u)^n.
\ee
The existence of such subsolutions provides the existence of a genuine solution, as established in \cite{SW,FLM} (see also \cite{Sze,CJY,PT} for extensions and generalizations).

\medskip
(e) There is another flow superficially similar to (\ref{AF}), but which can be considered for any compact complex manifold $X$, and not just manifolds which admit a nowhere holomorphic $(n,0)$-form $\Omega$,
\bea
\label{AF1}
\p_t\omega^{n-1}=i\p\bar\p\,\o^{n-2},
\eea
with initial data $\o_0$ satisfying $d\o_0^{n-1}=0$. 
A similar computation to the one in the proof of Theorem \ref{main4},
using Lemma 4 in \cite{PPZ5}, shows that the flow (\ref{AF1}) can be expressed as
\bea\label{af1-local}
\p_t g_{\bar{k} j} = -{1 \over (n-1) }\nabla^m T_{\bar k j m}+{1 \over 2 (n-1)} \{-g^{q \bar{p}} g^{s \bar{r}} T_{\bar k q s} \bar{T}_{j\bar p\bar r} + {|T|^2 \over n-1} g_{\bar k j}\}
\eea
for $n\geq 4$. From this formula, it is easy to deduce case (i) in Lemma \ref{cf+ak}, namely that the stationary points of the flow (\ref{AF1}) are K\"ahler metrics: it suffices to set (\ref{af1-local}) to zero and to take the trace in order to obtain
$|T|^2=0$. However, the flow (\ref{AF1}) may be hard to use, because it is not parabolic, and its stationary set may be too large, as it contains all K\"ahler metrics.

\medskip
(f) The flow (\ref{AF}) can be viewed as a K\"ahler analogue of the inverse Gauss curvature flow. Indeed, if we consider a flow of a strictly closed convex hypersurface $M_t$ in ${\bf R}^n$ by the inverse of its Gauss curvature, then it can be expressed \cite{Tsai} as
\bea
\p_t u
= {\det (u g_{ij}+ \na_i \na_j u ) \over \det g_{ij}}, \ \ \ u(x,0)=u_0(x)>0,
\eea
where $u$ is the support function $u: {\bf S}^n \times [0,T)$ defined by $u(N,t) = \langle P, N \rangle$ where $P \in M_t$ is the point on $M_t$ with normal $N$, and $({\bf S}^n,g_{ij})$ is the standard sphere. The right hand side of this flow exhibits the determinant of the Hessian of the unknown $u$, just as the right hand side of the equation (\ref{MA-flow}).

\begin{appendix}
\section{Appendix}
\setcounter{equation}{0}

In this Appendix, we group together our conventions for differential forms and several identities needed for the proof of Theorem \ref{main1} and Theorem \ref{main4}.

\subsection{Components of a differential form}
Let $\varphi$ be a $(p,q)$-form on the manifold $X$. We define its components
$\varphi_{\bar k_1\cdots\bar k_q j_1\cdots j_p}$ by
\bea
\varphi=
{1\over p!q!}
\sum \varphi_{\bar k_1\cdots\bar k_q j_1\cdots j_p}\,
dz^{j_p}\wedge\cdots\wedge dz^{j_1}\wedge
d\bar z^{k_q}\wedge\cdots\wedge d\bar z^{k_1}.
\eea

\subsection{Contraction identities}
We note a few basic identities for contracting differential forms of degree $(p,p)$. Let $\omega = i g_{\bar{k} j} dz^j \wedge d \bar{z}^k$ be a Hermitian metric. For a $(1,1)$ form $\alpha$
\be
\alpha = \alpha_{\bar{p} q} \, dz^q \wedge d \bar{z}^p,
\ee
we have
\be \label{1-1-contract}
\alpha \wedge {\omega^{n-1} \over (n-1)!} = -i (g^{j \bar{k}} \alpha_{\bar{k} j}) \, {\omega^n \over n!}.
\ee
Next, for a $(2,2)$ form $\Phi$ with components
\be
\Phi = {1 \over 4} \Phi_{\bar{p} q \bar{r} s} \, dz^s \wedge d \bar{z}^r \wedge dz^q \wedge d \bar{z}^p,
\ee
we have
\be \label{2-2-contract}
\Phi \wedge {\omega^{n-2} \over (n-2)!} = - {1 \over 2} \bigg\{ g^{j \bar{k}} g^{\ell \bar{m}} \Phi_{\bar{k} j \bar{m} \ell} \bigg\} \, {\omega^n \over n!}.
\ee
For a $(3,3)$ form $\Psi$ with components
\be
\Psi = {1 \over 36} \Psi_{\bar{k} j \bar{p} q \bar{r} s} \, dz^s \wedge d \bar{z}^r \wedge dz^q \wedge d \bar{z}^p \wedge dz^j \wedge d \bar{z}^k,
\ee
we have
\be \label{3-3-contract}
\Psi \wedge {\omega^{n-3} \over (n-3)!} = {i \over 6} \bigg\{ g^{j \bar{k}} g^{q \bar{p}} g^{s \bar{r}} \Psi_{\bar{k} j \bar{p} q \bar{r} s} \bigg\} \, {\omega^n \over n!}.
\ee

\subsection{Computing $T \wedge \bar{T}$}
Next, let $T$ be $(2,1)$ form. 
\be
T = {1 \over 2} T_{\bar{k} sj} dz^j \wedge dz^s \wedge d \bar{z}^k. 
\ee
Then
\be
T \wedge \bar{T} = {1 \over 4} T_{\bar{k} sj} \bar{T}_{q \bar{p} \bar{r}} \, d z^j \wedge d \bar{z}^k \wedge d z^s \wedge d \bar{z}^r \wedge dz^q \wedge d \bar{z}^p. 
\ee
Antisymmetrizing
\be
T \wedge \bar{T} = {1 \over (3!)^2} (T \wedge \bar{T})_{\bar{p} q \bar{r} s \bar{k} j} \, d z^j \wedge d \bar{z}^k \wedge d z^s \wedge d \bar{z}^r \wedge dz^q \wedge d \bar{z}^p. 
\ee
where
\bea\label{TwedgeT}
(T \wedge \bar{T})_{\bar{p} q \bar{r} s \bar{k} j} &=& \bigg\{ T_{\bar{k} sj} \bar{T}_{q \bar{p} \bar{r}} + T_{\bar{r} sj} \bar{T}_{q \bar{k} \bar{p}} +   T_{\bar{p} sj} \bar{T}_{q \bar{r} \bar{k}} +  T_{\bar{k} qs} \bar{T}_{j \bar{p} \bar{r}}+  T_{\bar{r} qs} \bar{T}_{j \bar{k} \bar{p}} \nonumber\\
&&+  T_{\bar{p} qs} \bar{T}_{j \bar{r} \bar{k}} +  T_{\bar{k} jq} \bar{T}_{s \bar{p} \bar{r}} +  T_{\bar{r} jq} \bar{T}_{s \bar{k} \bar{p}} + T_{\bar{p} jq} \bar{T}_{s \bar{r} \bar{k}} \bigg\}.
\eea
Let
\be
\tau = T_q dz^q, \ \ T_q = g^{j \bar{k}} T_{\bar{k} j q}, \ \ |\tau|^2 = g^{q \bar{p}} T_q \bar{T}_{\bar{p}}, \ \ |T|^2 = g^{q\bar p} g^{s\bar r} g^{j \bar{k}} T_{\bar p s j} \bar{T}_{q\bar r \bar k}.
\ee
Then
\bea\label{tr-tr-TT}
g^{q \bar{p}} g^{s \bar{r}} (T \wedge \bar{T})_{\bar{r} s \bar{p} q \bar{k} j} 
= g^{q\bar p} g^{s\bar r} ( 2T_{\bar p s j} \bar{T}_{q\bar r \bar k} + T_{\bar k q s} \bar{T}_{j\bar p\bar r}) 
-2 g^{s\bar r} (T_{\bar k j s} \bar{T}_{\bar r} + T_s \bar{T}_{j\bar k \bar r}) - 2 T_j \bar{T}_{\bar k},
\eea
and
\be\label{tr-tr-tr-TT}
g^{j \bar{k}} g^{q \bar{p}} g^{s \bar{r}} (T \wedge \bar{T})_{\bar{p} q \bar{r} s \bar{k} j} = 3 \{ |T|^2 - 2 |\tau|^2 \}.
\ee
Applying formula (\ref{3-3-contract}), we obtain
\be \label{T-Tbar-contract}
T \wedge \bar{T} \wedge \omega^{n-3} = {i \over 2} {1 \over n(n-1)(n-2)} \{ |T|^2 - 2 |\tau|^2 \} \, \omega^n .
\ee

\subsection{Computing $i \ddb \omega$}
We have
\be
i \ddb \omega = {1 \over 2^2} (i \ddb \omega)_{\bar{k} j \bar{\ell} m} \, dz^m \wedge d \bar{z}^\ell \wedge dz^j \wedge d \bar{z}^k,
\ee
with 
\be
(i \ddb \omega)_{\bar{k} j \bar{\ell} m} = \p_{\bar{\ell}} \p_j g_{\bar{k} m} - \p_{\bar{\ell}} \p_m g_{\bar{k} j} - \p_{\bar{k}} \p_j g_{\bar{\ell} m} + \p_{\bar{k}} \p_m g_{\bar{\ell} j}.
\ee
Using the definition of curvature,
\be
(i \ddb \omega)_{\bar{k} j \bar{\ell} m} = R_{\bar{k} j \bar{\ell} m} - R_{\bar{k} m \bar{\ell} j} + R_{\bar{\ell} m \bar{k} j} - R_{\bar{\ell} j \bar{k} m} - g^{s \bar{r}} T_{\bar{r} mj} \bar{T}_{s \bar{\ell} \bar{k}}.
\ee
Therefore
\be \label{tr-iddb-omega}
g^{j \bar{k}} (i \ddb \omega)_{\bar{k} j \bar{m} \ell} = \tilde{R}_{\bar{\ell} m} - R''_{\bar{\ell} m} + R_{\bar{\ell} m} - R'_{\bar{\ell} m} - g^{j \bar{k}} g^{s \bar{r}} T_{\bar{r} mj} \bar{T}_{s \bar{\ell} \bar{k}}. 
\ee

\end{appendix}

\bigskip
Department of Mathematics, Columbia University, New York, NY 10027, USA

\smallskip

phong@math.columbia.edu

\bigskip
Department of Mathematics, Columbia University, New York, NY 10027, USA

\smallskip
picard@math.columbia.edu

\bigskip
Department of Mathematics, University of California, Irvine, CA 92697, USA

\smallskip
xiangwen@math.uci.edu

\end{document}